\documentclass[reqno, 12pt]{amsart}
\usepackage{amssymb,amsmath,amsfonts}
\usepackage[margin=1in]{geometry}
\usepackage{dsfont}
\usepackage{color}
\usepackage{graphicx,graphics}
\usepackage{latexsym}
\usepackage[dvips]{epsfig}
\usepackage{mathtools,mathrsfs}
\usepackage{stmaryrd,cite}
\usepackage{cases}
\usepackage{enumerate}
\usepackage[usenames,dvipsnames,svgnames]{xcolor}
\usepackage{todonotes}
\usepackage{marginnote}
\usepackage{listings}
\usepackage{cleveref}
\allowdisplaybreaks
\definecolor{refkey}{rgb}{1,0,0.5}
\definecolor{labelkey}{rgb}{0,0.4,1}
\numberwithin{equation}{section}
\newtheorem{thm}{Theorem}[section]

\newtheorem{lem}[thm]{Lemma}

\newtheorem{rmk}[thm]{Remark}

\newcommand{\ea}{\epsilon}

\newcommand{\al}{\alpha}

\newcommand{\da}{\delta}

\newcommand{\ga}{\gamma}

\newcommand{\Oa}{\Omega}

\newcommand{\iy}{\infty}
\newcommand{\pl}{\partial}

\newcommand{\lt}{\left}
\newcommand{\rt}{\right}
\newcommand{\be}{\begin{equation}}
\newcommand{\ee}{\end{equation}}
\newcommand{\bee}{\begin{equation*}}
\newcommand{\eee}{\end{equation*}}

\newcommand{\ef}{\eqref}
\newcommand{\f}{\frac}

\begin{document}
\title[] {Time-Asymptotics of Physical Vacuum Free Boundaries for Compressible Inviscid Flows with Damping}
\author{Huihui Zeng}
\maketitle

\begin{abstract}
In this paper, we prove the leading term of
time-asymptotics of the moving vacuum boundary for compressible inviscid flows with damping
to be that for Barenblatt self-similar solutions to the corresponding porous media equations obtained by simplifying momentum equations via Darcy's law plus the possible shift due to the movement of the center of mass, in the one-dimensional and three-dimensional spherically symmetric motions, respectively.
This gives a complete description of the large time asymptotic behavior of solutions to the corresponding vacuum free boundary problems.
The results obtained  in this paper are the first ones   concerning  the large time asymptotics of physical vacuum boundaries for compressible inviscid fluids,
to the best of our knowledge.
 \end{abstract}

%{\small\textbf{Keywords}: }
	
%{\small \textbf{2010 Mathematics Subject Classification numbers}: }

% \tableofcontents

\section{Introduction}
The aim of this paper is to give the precise  time-asymptotics of physical   vacuum free boundaries for the following problem of compressible  Euler equations with damping:
\begin{subequations}\label{m2.1} \begin{align}
& \pl_t \rho   + {\rm div}(\rho   u ) = 0 &  {\rm in}& \ \  \Omega(t), \label{2.1a}\\
 &  \pl_t  (\rho   u )   + {\rm div}(\rho   u \otimes   u )+\nabla_{  x} p(\rho) = -\rho {  u}  & {\rm in}& \ \ \Omega(t),\label{2.1b}\\
 &\rho>0 &{\rm in }  & \ \ \Omega(t),\label{2.1c}\\
 & \rho=0    &    {\rm on}& \  \ \Gamma(t)=\pl \Omega(t), \label{2.1d}\\
 &    \mathcal{V}(\Gamma(t))={  u}\cdot \mathcal{N}, & &\label{2.1e}\\
&(\rho,{ u})=(\rho_0, { u}_0) & {\rm on} & \ \   \Omega(0), \label{2.1f}
 \end{align} \end{subequations}
where $( x,t)\in \mathbb{R}^n \times [0,\iy) $,  $\rho $, ${ u} $, and $p$ denote, respectively, the space and time  variable, density, velocity and  pressure; $\Omega(t)\subset \mathbb{R}^n$, $\Gamma(t)$, $\mathcal{V}(\Gamma(t))$ and $ \mathcal{N}$ represent, respectively, the changing domain occupied by the gas at time $t$, moving vacuum boundary, normal velocity of $\Gamma(t)$, and exterior unit normal vector to $\Gamma(t)$.  We are concerned with the polytropic gas for which the equation of state is  given by
$$ p(\rho)=\rho^{\gamma}, \  \ {\rm where} \ \  \gamma>1 {~\rm is~the~adiabatic~exponent}. $$
Let $c(\rho)=\sqrt{ p'(\rho)}$ be the sound speed, the condition
\be\label{physical vacuum} -\infty<\nabla_\mathcal{N}\lt(c^2(\rho)\rt)<0  \  \ {\rm on} \ \  \Gamma(t) \ee
defines {\em a  physical vacuum boundary}  (cf. \cite{7,10',16',23,24,25}), which is also called {\em a vacuum boundary with physical singularity} in contrast to the case that $ \nabla_\mathcal{N}\lt(c^2(\rho)\rt)=0$ on  $\Gamma(t)$.
The role of physical vacuum singularity plays is to push vacuum boundaries. Indeed, one has that $D_t u \cdot \mathcal{N} = - (\ga-1)^{-1} \nabla_\mathcal{N}\lt(c^2(\rho)\rt) - u\cdot \mathcal{N}$ on $\Gamma(t)$, by restricting the momentum equation \ef{2.1b} on $\Gamma(t)$, where $D_t u= (\pl_t + u\cdot \nabla_x)u$ is the acceleration of $\Gamma(t)$, and the term  $- (\ga-1)^{-1} \nabla_\mathcal{N}\lt(c^2(\rho)\rt)>0$ represents the pressure effect to accelerate vacuum boundaries.
In order to capture this physical singularity, the initial density is supposed to satisfy
\begin{subequations}\label{initial density}
\begin{align}
& \rho_0>0 \ \ {\rm in} \ \ \Omega(0), \ \  \rho_0=0 \ \ {\rm on} \ \ \Gamma(0),   \ \ \int_{\Omega(0)} \rho_0(x) dx =M,
\\
&  -\infty<\nabla_\mathcal{N}\lt(c^2(\rho_0)\rt)<0  \  \ {\rm on} \ \  \Gamma(0), \end{align}
\end{subequations}
where $M\in (0, \iy)$ is the initial total mass.

The compressible Euler equations of isentropic flows  with damping, \ef{2.1a}-\ef{2.1b}, is closely related to the
porous media equation (cf. \cite{HL, HMP, HPW, 23, LZ, HZ}):
\begin{equation}\label{pm}
\pl_t \rho  =\Delta p(\rho),
\end{equation}
when \ef{2.1b} is simplified to Darcy's law:
\begin{equation}\label{darcy}
\nabla_{ x} p(\rho)=- \rho { u}.
\end{equation}
(The equivalence can  be seen formally by the rescaling ${ x}'=\ea { x}, t'=\ea^2 t, { u}'= { u}/ \ea$.)
Basic understanding of the solution to \ef{pm} with finite mass is provided by Barenblatt (cf. \cite{ba}), which is given by
\begin{align}\label{Feb25-1}
\bar \rho({ x}, t)=
(1+  t)^{-\frac{n}{n\ga-n+2}}\left(A_n- B_n(1+ t)^{-\frac{2}{n\ga-n+2}}|{ x}|^2\right)^{\frac{1}{\gamma-1}},
\end{align}
defined in
\begin{align}\label{Feb25-2}
\bar \Omega(t)= \lt\{x\in\mathbb{R}^n:\  |x| < A_n^{\frac{1}{2}} B_n^{-{\frac{1}{2}}}(1+t)^{\frac{1}{n\ga-n+2}} \rt\},
\end{align}
where $A_n$ and $B_n$ are positive constants determined by $\ga$ and $M$. It holds that for $t\ge 0$,
$$\bar\rho> 0 \ \ {\rm in} \ \  \bar \Omega(t), \ \  \bar\rho=0 \ \ {\rm on} \ \ \bar \Gamma(t)= \pl \bar \Omega(t), \ \ {\rm and} \ \   \int_{\bar \Omega(t)}  \bar\rho(t,{ x} )d{ x}  =M . $$
We define the corresponding  Barenblatt velocity $\bar { u}$ by
\begin{align}
\bar { u}(x,t)=- \f{ \nabla_{ x} p(\bar \rho) }{\bar\rho}
=\frac{1}{n\ga-n+2}\f{  x }{1+t}
 \ \ {\rm in} \ \  \bar \Omega (t). \label{Feb26}
\end{align}
So,  \ef{pm}-\ef{darcy}   are solved by $(\bar\rho, \bar { u})$ defined in the region $\bar \Omega(t)$. When $\gamma$ is fixed, there is only one
parameter, total mass $M$, for the Barenblatt self-similar solution.
It is  assumed that the initial total mass of problem \ef{m2.1}  is the same as that of the Barenblatt solution:
\be\label{Feb25-4}
\int_{\Omega(0)} \rho_0(x) dx =M=\int_{\bar \Omega(t)}  \bar\rho(t,{ x} )d{ x}.
\ee

Apparently,  the physical vacuum condition \ef{physical vacuum} is satisfied by the Barenblatt solution. This is the major motivation to study problem \ef{m2.1} with the initial condition \eqref{initial density}.  To this end,   a family of particular solutions to problem \ef{m2.1} was constructed in \cite{23} by use of the following spherically-symmetric ansatz:
\begin{subequations}\label{liuexplicitsolution}\begin{align}
&\Omega(t) = \lt\{x\in\mathbb{R}^n:\  |x| < \sqrt {e(t)/b(t)}  \rt\} ,  \\
&  c^2({x}, t)={e(t)-b(t)|x|^2},  \ \
  { u}({ x}, t)= a(t) x.
  \end{align}\end{subequations}
In \cite{23}, a  system of ordinary differential equations for $(e,b,a)(t) $ was derived  with $e(t), b(t)>0$ for $t\ge 0$, and it was shown that this family of particular solutions  is  time-asymptotically equivalent to the Barenblatt self-similar solution  with the same total mass.
By the same ansatz as \ef{liuexplicitsolution}, the Barenblatt solution  of \ef{pm}-\ef{darcy} can be obtained as $
 \bar c^2({ x}, t)=\bar e(t)-\bar b(t)r^2$ and  $  { u}({ x}, t)= \bar a(t) {{ x}}  $, and
it was shown in \cite{23}  that
  $$
  (a,\ b,\ e)(t)=(\bar a, \ \bar b, \ \bar e)(t)+ O(1)(1+t)^{-1}{\ln (1+t)} \ \ {\rm as}\  \ t\to\infty.$$
It has been an important open question since the  construction of particular solutions to \ef{m2.1} in \cite{23}, that,  for the general initial data, whether there is still a large time existence theory for  problem \ef{m2.1} capturing the physical vacuum singular behavior \ef{physical vacuum}, and if the  time-asymptotic equivalence still holds for the solution to \ef{m2.1} and the corresponding Barenblatt self-similar solution with the same total mass.
This issue was addressed,  respectively,   in the one-dimensional case (cf. \cite{LZ}) and three-dimensional spherically symmetric case (cf. \cite{HZ}), by proving the
global existence of smooth solutions to  problem  \ef{m2.1} and the convergence  of the density and velocity to those of Barenblatt solutions.
The existence of almost global solution was recently proved in \cite{Z3} for the general three-dimensional perturbations without any symmetry assumptions.
It was also proved in \cite{LZ} and \cite{HZ}, respectively, that the expanding rate of the vacuum boundary $\Gamma(t)$ is exactly the same as that of the corresponding Barenblatt self-similar solution, of the order of $(1+t)^{{1}/({n\ga-n+2})}$.  However, the large time behavior of $\Gamma(t)$ was not given in either
\cite{LZ} or \cite{HZ}.
The purpose of this paper is to   fill  this gap by proving the precise time-asymptotics of the moving vacuum boundary  to give a complete description of the large time asymptotic behavior of solutions to  problem \ef{m2.1} for the one-dimensional and three-dimensional spherically symmetric motions, respectively.

The main results of this paper  concern the precise time-asymptotics of the moving vacuum boundary $\Gamma(t)$.
As $t\to \infty$, the leading term of  $\Gamma(t)$ is  the Barenblatt vacuum boundary $ \bar \Gamma(t)$ plus a shift determined by the initial center of mass and average velocity. More explicitly, we obtain for any constant $q>0$,
\bee
\Gamma(t)=\bar \Gamma(t)+ \vartheta_0+O(1)\lt(\lt((\ln(1+t))^q+1\rt)^{-1}\sqrt{\mathcal{E}(0)}+ (1+t)^{ -1+\frac{1}{n\ga-n+2}}\ln(1+t)\rt)
\eee
where $O(1)$ is a bounded quantity depending only on $\ga$,  $M$ and $q$;
$\mathcal{E}(0)$ is the norm for the initial perturbation to the Barenblatt solution whose definition will be given in \ef{d3.29-1} and \ef{d3.29-2};
 $\vartheta_0$ is the sum of the initial center of mass and average velocity, given by
\be\label{d3.12-1}\vartheta_0=\frac{1}{M} \int_{\Oa(0)} x\rho_0(x)dx +\frac{1}{M}  \int_{\Oa(0)} \rho_0(x) u_0(x) dx.
\ee
It should be noted that $\vartheta_0=0$ for the spherically symmetric initial data. We also obtain better convergence rates in this paper of the density and velocity for problem \ef{m2.1} to the corresponding ones for Barenblatt solutions than those given in
 \cite{LZ,HZ}.

The physical vacuum phenomena that the sound speed $c(\rho)$ is $C^{ {1}/{2}}$-H$\ddot{\rm o}$lder continuous across  vacuum boundaries appear in several important physical situations such as gaseous stars (cf. \cite{ HaJa2, HaJa3, LXZ}), shallow waters, besides gas flows through porous media modelled by \ef{m2.1}.
The mathematical challenge of  physical vacuum singularities lies in the fact that  the characteristic speeds become singular with infinite spatial derivatives
at  vacuum boundaries, and  the standard approach of symmetric hyperbolic systems (cf. \cite{Friedrichs, Kato,17}) do not apply.
Early study of well-posedness of smooth solutions with sound speed $c(\rho)$ smoother than $C^{ {1}/{2}}$-H$\ddot{\rm o}$lder continuous  near vacuum states for compressible inviscid fluids can be found in \cite{chemin1,chemin2,24,25,MUK,Makino,39}.
The strong  degeneracy of equations near vacuum states makes it extremely difficult to obtain regularity estimates near physical vacuum boundaries.
Indeed, due to its physical significance and challenge in mathematical analyses, physical vacuum problems of compressible Euler and related equations have received much attentions.
The local-in-time well-posedness theory is established  in \cite{16,10, 7, 10', 16'}. (See also \cite{zhenlei1, LXZ,serre} for related works on the local theory.)
Understanding the large time stability of solutions with scaling invariance such as Barenblatt self-similar solutions and affine motions is of fundamental importance in the study of physical vacuum singularities,  which poses great challenges due to extreme difficulties in obtaining the uniform-in-time  higher order regularity of solutions near vacuum boundaries.
In this direction, besides the results in \cite{LZ, HZ, Z3} mentioned earlier,  there have been global-in-time results for expanding solutions to compressible Euler equations (cf. \cite{HaJa1,ShSi}) when the initial data are small perturbations of affine motions (cf. \cite{sideris1,sideris2}), and to compressible Euler-Poisson equations of gaseous stars (cf. \cite{HaJa2, HaJa3}).
It should be noted that the expanding rates of vacuum boundaries of expanding solutions are linear, $O(1+t)$,
but those of Barenblatt solutions, which are the background approximate solutions  for \ef{m2.1} in large time, are sub-linear, $O((1+t)^{{1}/({n(\ga-1)+2})}) $.
This slower expanding rate of  vacuum boundaries induces the slower decay of various quantities.

For the Cauchy problem of the one-dimensional compressible Euler equations with damping,  the $L^p$-convergence of
$L^{\infty}$-weak solutions to Barenblatt solutions was given in \cite{HMP} with  $p=2$ if $1<\ga\le 2$ and $p=\ga$ if $\ga>2$ and in \cite{HPW} with $p=1$, respectively, using  entropy-type estimates for the solution itself without deriving  estimates for derivatives.
However, the vacuum boundaries separating gases and vacuum cannot be traced in the framework of $L^{\infty}$-weak solutions.

{\bf Notations}. we use $C(\beta)$ to denote  a certain positive constant
depending on quantity $\beta$. They are referred as universal and can change
from one inequality to another one.

\section{The one-dimensional motions}\label{section2}
In  this section, we consider the vacuum free boundary  problem \ef{m2.1} in the one-dimensional case. Since the changing domain $\Omega(t)$ is an interval in this situation, we write it as
$$\Oa(t)=\lt\{x\in\mathbb{R}:\  x_-(t)\le x \le x_+(t)\rt\}, \ \  t\ge 0 .$$
As in \cite{LZ},
we transform problem \ef{m2.1} into Lagrangian variables to fix the boundary,  and make the initial interval of the Barenblatt solution, $\bar \Oa(0)$, as the reference interval. For simplicity of presentation, set
$$ \mathcal{I} = \bar \Oa (0) =\lt(-\sqrt{A_1 /B_1 }, \ \sqrt{A_1 /B_1} \rt).  $$
For $x\in \mathcal{I}$, we define the  Lagrangian variable $\eta(x, t)$  by
\be\label{Feb28-5}
\eta_t(x, t)= u(\eta(x, t), t) \ \ {\rm for} \  \ t>0, \  \ {\rm and} \ \  \eta(x, 0)=\eta_0(x).
\ee
Then, it follows from the conservation of mass,  \ef{2.1a}, that
$\pl_t \lt( \rho(\eta(x, t), t)\eta_x(x, t) \rt)=0$, which implies
 $\rho(\eta(x, t), t)\eta_x(x, t) =\rho_0(\eta_0(x))\eta_{0}'(x)
$.
If we choose $\eta_0$ such that $\rho_0(\eta_0(x))\eta_{0}'(x) =    \bar\rho(x,0) $, the initial density of the Barenblatt solution, then we have for $x\in \mathcal{I}$,
\be\label{Feb25-3}
\rho(\eta(x, t), t)=\frac{\bar \rho_0(x)}{\eta_x(x, t)}  , \ \ {\rm where} \ \ \bar \rho_0(x)  = \bar\rho(x,0)=(A_1-B_1x^2)^{\frac{1}{\ga-1}}.
\ee
Such an $\eta_0$ exists, for example, we may define $\eta_0: \  \mathcal{I} \to \Oa(0)$ by
$$\int_{x_-(0)}^{\eta_0(x)}\rho_0(y)dy=\int_{-\sqrt{A_1/ B_1 }}^{x}\bar\rho_0(y)dy  \ \ {\rm for} \  \ x \in \mathcal{I} .$$
Due to \ef{initial density}, \ef{Feb25-1} and \ef{Feb25-4}, $\eta_0$ defined above is a well-defined diffeomorphism.
So, the vacuum free boundary problem \ef{m2.1} is reduced to the following initial-boundary value problem with the fixed
boundary:
\begin{subequations}\label{equ-1d}\begin{align}
& \bar\rho_0 \eta_{tt} + \bar\rho_0 \eta_{t}+ \lt( \bar\rho_0^\ga \eta_x^{-\ga}  \rt)_x= 0  \ \  &{\rm in} \  \   \mathcal{I}\times (0, \iy), \label{equa}\\
&(\eta, \eta_t)= \lt( \eta_0,  u_0(\eta_0)\rt)  \ \  &   {\rm on}  \  \  \mathcal{I} \times\{t=0\}.
\end{align}\end{subequations}
In the setting, the moving vacuum boundaries for problem \ef{m2.1} are given by
\be\label{vbs-1d}  x_{\pm}(t)=\eta\lt(\pm\sqrt{A_1/ B_1},t \rt) \ \  {\rm for} \ \  t\ge 0.\ee

\subsection{Main results}
Define the Lagrangian variable $\bar\eta(x, t)$ for the Barenblatt flow in $\bar {\mathcal{I}}$ by
$$
 \bar\eta_t (x, t)= \bar u(\bar\eta(x, t), t)=\f{ \bar\eta(x,t)}{(\ga+1)(1+  t)} \ \ {\rm for} \  \ t>0 \  \ {\rm and} \ \  \bar\eta(x, 0)=x,
$$
then
\be\label{212}
\bar\eta(x,t)=x(1+ t)^{\f{1}{\ga+1}} \ \ {\rm for} \ \  (x,t) \in \bar {\mathcal{I}}\times [0,\iy).
\ee
Since $\bar\eta$ does not solve $\ef{equa}$ exactly, a correction $h(t)$ was introduced in \cite{LZ} which solves the following initial   value problem of ordinary differential equations:
\begin{align*}
  &h_{tt} + h_t   -  (\ga+1)^{-1} {(\bar\eta_x+h)^{-\ga}} +  \bar\eta_{xtt}   + \bar\eta_{xt} =0 ,  \ \ t>0, \\
  &h(t=0)=h_t(t=0)=0.
\end{align*}
The ansatz is then set   to be
$
\tilde{\eta}(x,t) =\bar\eta(x,t)+x h(t),
$
so that
\be\label{equeta-1d} \begin{split}
  \bar\rho_0   \tilde\eta_{tt}   +    \bar\rho_0\tilde\eta_t + \lt(\bar \rho_0^\ga \tilde \eta_x^{-\ga}  \rt)_x= 0  \ \  {\rm in} \  \   \mathcal{I}\times (0, \iy).
\end{split}
\ee
It was shown in \cite{LZ} that $\tilde\eta$ behaves similarly to $\bar \eta$, see \ef{decay-1d} and \ef{h-1d} for  details.  We write equation \ef{equa} in a perturbation form around the Barenblatt solution.
Let $w(x,t)=\eta(x,t)- \tilde\eta(x,t)$,
then it follows from \ef{equa} and \ef{equeta-1d} that
\be\label{eluerpert}\begin{split}
\bar\rho_0w_{tt}   +
   \bar\rho_0w_t + \lt[ \bar\rho_0^\ga \lt(  (\tilde\eta_x+w_x)^{-\ga} -   {\tilde\eta_x}^{-\ga}  \rt) \rt]_x =0.
\end{split}\ee
We denote  $\alpha=(\ga-1)^{-1}$, and set that for $j=0,\cdots, 4+ [\alpha]$ and  $i=1,\cdots, 4+[\alpha]-j$,
\bee\label{}\begin{split}
&\mathcal{ E}_{j}(t)   =   (1+  t)^{2j} \int_\mathcal{I}    \lt[\bar\rho_0\lt(\pl_t^j w\rt)^2 + \bar\rho_0^\ga \lt(\pl_t^j  w_x \rt)^2 + (1+ t) \bar\rho_0   \lt(\pl_t^{j} w_t\rt)^2 \rt] (x, t)  dx  ,  \\
&\mathcal{ E}_{j, i}(t)  =   (1+  t)^{2j}  \int_\mathcal{I}  \lt[\bar\rho_0^{1+(i-1)(\ga-1) }  \lt(\pl_t^j \pl_x^i w\rt)^2 + \bar\rho_0^{\ga+ i(\ga-1) } \lt(\pl_t^j  \pl_x^{i}w_x  \rt)^2\rt] (x, t) dx .
\end{split}\eee
The higher-order norm is defined by
\be\label{d3.29-1}
  \mathcal{E}(t)  =  \sum_{0\le j\le 4+[\alpha]} \lt(\mathcal{ E}_{j}(t) + \sum_{1\le i\le 4+[\alpha]-j} \mathcal{ E}_{j, i}(t)  \rt).
\ee

The main results in this section are stated as follows.

\begin{thm}\label{1d-mainthm} There exists a  constant $\ea_0 >0$ such that if
$\mathcal{E}(0)\le \ea_0^2$,
then   problem \eqref{equ-1d}  admits a global unique smooth solution  in $\mathcal{I}\times[0, \iy)$ satisfying for all  $t\ge 0$,
\begin{subequations}\label{Feb28-1}\begin{align}
&\mathcal{E}(t)\le C(\ga, M)\mathcal{E}(0),\\
 &\lt(\lt(\ln(1+t)\rt)^q +1\rt) \lt( \widetilde{\mathcal{E}}(t)  +  \int_\mathcal{I} \bar\rho_0 (w-\vartheta_0)^2 dx \rt)
\le C(\ga, M,q)\mathcal{E}(0)\label{d3.12}
\end{align}\end{subequations}
for any fixed constant $q>0$,
where $\widetilde{\mathcal{E}}(t)=\mathcal{E}(t)-\int_\mathcal{I} \bar\rho_0 w^2 dx$, and
$\vartheta_0$ is the sum of the initial center of mass and  average velocity, given by \ef{d3.12-1}.

\end{thm}

With  Theorem \ref{1d-mainthm},  we have the following theorem for solutions to the original vacuum free boundary problem
\ef{m2.1} concerning the time-asymptotics of vacuum boundaries $x_{\pm} (t)$ and better convergence rates of the density $\rho$ and velocity $u$
than those in \cite{LZ}.

{\begin{thm}\label{1d-mainthm1} There exists a  constant $\ea_0 >0$ such that if
$\mathcal{E}(0)\le \ea_0^2$,
then   problem
\ef{m2.1} with $n=1$  admits a global unique smooth solution $\lt(\rho, u,  \Oa(t)\rt)$ for $t\in[0,\iy)$. Let $q>0$ be any fixed constant, and set
$$\mathscr{G}(t)=\lt((\ln(1+t))^q+1\rt)^{-1}\sqrt{\mathcal{E}(0)}+ (1+t)^{-\frac{\ga}{\ga+1}}\ln(1+t),
$$
then it holds that for all $x \in \mathcal{I}$ and $t\ge 0 $,
\begin{subequations}\label{20.3.28}\begin{align}
& \lt|\rho\lt(\eta(x, t),t\rt)-\bar\rho\lt(\bar\eta(x, t), t\rt)\rt|
 \le   C(\ga, M,q)\lt(A_1-B_1x^2\rt)^{\frac{1}{\ga-1}}(1+t)^{-\frac{2}{\ga+1}}  \mathscr{G}(t) , \label{1'-1d}\\
   &  \lt|u\lt(\eta(x, t),t\rt)-\bar u\lt(\bar\eta(x, t), t\rt)\rt|  \le C(\ga, M,q)(1+t)^{-1} \mathscr{G}(t); \label{2'-1d}
\end{align}\end{subequations}
and for all $t\ge 0$,
\begin{subequations}\begin{align}
& \lt|x_{\pm}(t)-\bar x_\pm( t) - \vartheta_0 \rt|  \le   C(\ga,M,q) \mathscr{G}(t), \label{Feb28-3}\\
& \lt|\frac{d^k x_{\pm}(t)}{dt^k}\rt|\le C(\ga, M)(1+t)^{\frac{1}{\ga+1}-k} , \ \   k=1, 2, 3 , \label{4'-1d}
\end{align}\end{subequations}
where $\bar x_\pm( t)=\pm \sqrt{A_1/ B_1} (1+ t)^{ {1}/({\gamma+1})} $
are the vacuum  boundaries of the Barenblatt solution, and $\vartheta_0$ is defined by \ef{d3.12-1}.
\end{thm}

\begin{rmk} In \cite{LZ}, it was proved
that    the estimates in \ef{20.3.28}   hold with
$$\mathscr{G}(t)=\sqrt{\mathcal{E}(0)}+ (1+t)^{-\frac{\ga}{\ga+1}}\ln(1+t).
$$
So,  the estimates in  \ef{20.3.28}  improve  the convergence rates of the density and velocity obtained in \cite{LZ}. Compared with \ef{Feb28-3}, it was shown in \cite{LZ} that
for all $t\ge 0 $,
\begin{align}
 \lt|x_{\pm}(t)-\bar x_\pm( t)  \rt|  \le   C(\ga,M)\lt(\sqrt{\mathcal{E}(0)}+ (1+t)^{-\frac{\ga}{\ga+1}}\ln(1+t)\rt), \label{20.3.28-c}
\end{align}
which means that $x_\pm(t)$ and $\bar x_\pm( t)$, vacuum boundaries of problem \ef{m2.1} and Barenblatt solutions, have the same expanding rates.
In addition to this,  estimate \ef{Feb28-3} gives the precise time-asymptotics of $x_\pm(t)$.
\end{rmk}

\subsection{Proof of Theorem \ref{1d-mainthm}}
It follows from Theorem 2.1 in \cite{LZ} that
there exists a small constant $\epsilon_0>0$ such that if
$\mathcal{E}(0)\le \epsilon_0^2$,
then    problem \eqref{equ-1d}  admits a global unique smooth solution $\eta(x,t)=w(x,t)+ \tilde\eta(x,t)$ in $\mathcal{I}\times[0, \iy)$ satisfying for all $t\ge 0$,
\bee\label{date2.3-1}
\mathcal{E}(t)\le C(\ga,M)\mathcal{E}(0)
\eee
and
\begin{align}\label{20lem31est}
  & \sum_{0\le j\le 3} (1+  t)^{2j} \lt\|\pl_t^j w(\cdot,t)\rt\|_{L^\iy(\mathcal{I})}^2   +  \sum_{j=0,1 } (1+  t)^{2j} \lt\|\pl_t^j w_x(\cdot,t)\rt\|_{L^\iy(\mathcal{I})}^2
   \notag \\
   & +
  \sum_{
  i+j\le 4+[\alpha],\ \  2i+j \ge 4 } (1+  t)^{2j}\lt\|  \bar\rho_0^{\f{(\ga-1)(2i+j-3)}{2}} \pl_t^j \pl_x^i w(\cdot,t)\rt\|_{L^\iy(\mathcal{I})}^2 \notag\\
  & \le C(\ga, M) \mathcal{E}(0) \le C(\ga, M) \epsilon_0^2.
\end{align}
This section is devoted to proving estimate \ef{d3.12}. For this purpose, we rewrite \ef{eluerpert} as
\begin{align}\label{12.30}
\bar\rho_0w_{tt}   +
   \bar\rho_0w_t - \ga {\tilde\eta_x}^{-\ga-1}  \lt( \bar\rho_0^\ga  w_x \rt)_x + \lt( \bar\rho_0^\ga R \rt)_x  =0,
   \end{align}
where $R$ is equivalent to $(1+t)^{-1-1/({\ga+1})} w_x^2$, given by
\be\label{date2.10}
{R} =   (\tilde\eta_x+w_x)^{-\ga} -   {\tilde\eta_x}^{-\ga} + \ga {\tilde\eta_x}^{-\ga-1} w_x.
\ee

As a starting point, we prove the following lemma.
\begin{lem}\label{lem2020} We have for all $t\ge 0$,
\begin{align}\label{1.31}
\mathscr{E}_j(t)+\int_0^t \mathscr{D}_j(s)ds \le C(\ga, M)\sum_{0\le m\le j}\mathscr{E}_m(0), \ \  j=0,1,\cdots,4+\lt[\alpha\rt],
\end{align}
where
\begin{align*}
\mathscr{E}_j(t)=&\int_\mathcal{I} \bar\rho_0^{1+j(\ga-1)}|\pl_x^j w|^2 dx
+(\ln(1+t)+1)\int_\mathcal{I} \bar\rho_0^{\ga+j(\ga-1)}|\pl_x^j w_x|^2 dx
\\
&+(1+t)(\ln(1+t)+1)\int_\mathcal{I} \bar\rho_0^{1+j(\ga-1)}|\pl_x^j w_t|^2 dx,\\
\mathscr{D}_j(t)=& (1+t)^{-1}\int_\mathcal{I} \bar\rho_0^{\ga+j(\ga-1)}|\pl_x^j w_x|^2 dx
\\
&+(1+t)(\ln(1+t)+1)\int_\mathcal{I} \bar\rho_0^{1+j(\ga-1)}|\pl_x^j w_t|^2 dx.
\end{align*}
\end{lem}
Once this lemma is proved, we will obtain the following estimate:
\begin{align}
 &(\ln(1+t)+1) \sum_{1\le i\le 4+[\alpha]} \mathcal{ E}_{0, i}(t) \le C(\ga,M)\mathcal{E}(0),
 \label{date2.3}
\end{align}
with the aid of the Hardy inequality \ef{hardy2019-1d}. It should be noted that estimate \ef{date2.3} can give the  time-asymptotics of the vacuum free boundaries as follows:
\begin{align}
 \lt|x_{\pm}(t)-\bar x_\pm( t) -\vartheta_0  \rt|  \le &  C(\ga,M) \lt(\ln(1+t)+1\rt)^{-1/{2}}\sqrt{ \mathcal{E}(0)} \notag \\
& +  C(\ga,M) (1+t)^{-{\ga}/({\ga+1})}\ln(1+t), \label{d3.24}
\end{align}
where  $\bar x_\pm( t)=\pm \sqrt{A_1/ B_1} (1+ t)^{ {1}/({\gamma+1})} $
are the vacuum  boundaries of the Barenblatt solution, and $\vartheta_0$ is the constant defined by \ef{d3.12-1}. Indeed, the derivation of \ef{d3.24} is the same as that of estimate \ef{Feb28-3}.
We can improve estimate \ef{date2.3} by use of the following lemma.

\begin{lem}\label{lem-d3.8}
Let $q$ be any fixed positive integer, then we have for all $t\ge 0$,
  \begin{align}
\lt(\lt( \ln(1+t) \rt)^q +1 \rt)\mathcal{E}_j(t; w-\vartheta_0) +  \int_0^t \lt( \lt( \ln(1+s) \rt)^q +1 \rt) \mathcal{D}_{j}(s) ds  \notag \\
  \le C(\ga, M,q)   \mathcal{E} (0), \ \  \ \  j=0 ,1 ,2, \cdots,  4 +[\alpha],  \label{new12-4}
\end{align}
where $\vartheta_0$ is the constant defined by \ef{d3.12-1},
\begin{align*}
\mathcal{ E}_{j}(t;w)   =  & (1+  t)^{2j} \int_\mathcal{I}    \lt[\bar\rho_0\lt(\pl_t^j w\rt)^2 + \bar\rho_0^\ga \lt(\pl_t^j  w_x \rt)^2  \rt] (x, t)  dx\\
&  + (1+  t)^{2j+1} \int_\mathcal{I}  \bar\rho_0   \lt(\pl_t^{j} w_t\rt)^2 (x, t)  dx,
\\
\mathcal{D}_{j}(t)   = & (1+ t)^{2j-1} \int_\mathcal{I}   \bar\rho_0^\ga  \lt|\pl_t^j w_x  \rt|^2  dx + (1+ t)^{2j+1} \int_\mathcal{I} \bar\rho_0   |\pl_t^j w_t|^2   dx.
\end{align*}
\end{lem}
Once this lemma is proved, estimate \ef{d3.12} will be a consequence of the fact that $ ( \ln(1+t)  )^q \le C(q)  ( ( \ln(1+t)  )^{[q]+1} +1 )$ for any positive $q\notin \mathbb{Z}$; the following elliptic estimates:
\begin{align}
   \mathcal{E}_{j, i}(t)\le C(\ga, M) \lt( \widetilde{\mathcal{E}}_0(t) + \sum_{1\le \iota \le i+j} \mathcal{E}_\iota(t) \rt) \label{20.el}
\end{align}
for $j\ge 0$, $ i\ge 1$, $  i+j\le 4 +[\alpha]$,
where $\widetilde{\mathcal{E}}_0(t)={\mathcal{E}}_0(t)-\int_{\mathcal{I}} \bar\rho_0 w^2 dx$; and
\begin{align*}
&\mathcal{E}_0(t; w-\vartheta_0)=\widetilde{\mathcal{E}}_0(t) + \int_{\mathcal{I}} \bar\rho_0(w-\vartheta_0)^2 dx, \\
&\mathcal{E}_j(t; w-\vartheta_0)=\mathcal{E}_j(t; w)= \mathcal{E}_j(t), \ \ 1\le j\le 4+[\alpha].
\end{align*}
Indeed, the elliptic estimate \ef{20.el} follows from
Proposition 3.1 in \cite{LZ}.

To prove Lemmas \ref{lem2020} and \ref{lem-d3.8},  we  need some  embedding estimates.
Let $\mathcal{U}$ be a bounded smooth domain in $\mathbb{R}^n$, and
 $d=d(x)=dist(x, \mathcal{U})$ be the distance function to the boundary.
For any  $a>0$ and nonnegative integer $b$, we define the  weighted Sobolev space  $H^{a, b}(\mathcal{U})$  by
$$ H^{a, b}(\mathcal{U})  = \lt\{   d^{a/2}F\in L^2(\mathcal{U}): \ \  \int_\mathcal{U}   d^a|\pl_x ^k F|^2dx<\infty, \ \  0\le k\le b\rt\}$$
  with the norm
 $ \|F\|^2_{H^{a, b}(\mathcal{U})}  = \sum_{k=0}^b \int_\mathcal{U}   d^a|\pl_x^k F|^2dx  $.
Then  for $b\ge  {a}/{2}$, we have the following {\it embedding of weighted Sobolev spaces} (cf. \cite{18'}):
 $ H^{a, b}(\mathcal{U})\hookrightarrow H^{b- {a}/{2}}(\mathcal{U})$
    with the estimate
  \be\label{wsv} \|F\|_{H^{b- {a}/{2}}(\mathcal{U})} \le C(a, b,\mathcal{U} ) \|F\|_{H^{a, b}(\mathcal{U})}. \ee
The following general version of the {\it Hardy inequality}, whose proof can be found  in \cite{18'}, will also be used in this paper. Let $k>1$ be a given real number and $F$ be a function satisfying
$
\int_0^{\da} x^k\lt(F^2 + F_x^2\rt) dx < \iy,
$
where $\da$ is a positive constant, then it holds that
\be\label{hardy2019}
\int_0^{\da} x^{k-2} F^2 dx \le C(\delta, k) \int_0^{\da} x^k \lt( F^2 + F_x^2 \rt)  dx.
\ee
This, together with  $\bar\rho_0^{\ga-1}=A-Bx^2$, implies
\be\label{hardy2019-1d}
\int_\mathcal{I}  \bar\rho_0^{(k-2)(\ga-1)}  F^2 dx \le C(\ga, M,k) \int_\mathcal{I}  \bar\rho_0^{k(\ga-1) } \lt( F^2 + F_x^2 \rt)  dx,
\ee
provided that the right hand side of \ef{hardy2019-1d} is finite.

\subsubsection{ Proof of Lemma \ref{lem2020}}
The proof consists of two steps. In {\em Step 1}, we show \ef{1.31} for $j=0$, the zeroth order estimate. Based on this, we use the mathematical induction to prove \ef{1.31} for $j\ge 1$, the higher order estimates in {\em Step 2}. First, we  recall  that for all  $t\ge 0$,
\begin{subequations}\label{decay-1d}\begin{align}
&\lt(1 +   t \rt)^{\frac{1}{\ga+1}} \le \tilde \eta_{x} (t)\le C(\ga) \lt(1 +   t \rt)^{\frac{1}{\ga+1}}, \ \  \ \  \tilde\eta_{xt}(t)\ge 0,  \label{decay-1d-a}\\
&\lt|\f{  d^k\tilde \eta_{x}(t)}{dt^k}\rt| \le C(\ga, k)\lt(1 +   t \rt)^{\frac{1}{\ga+1}-k},   \ \ k= 1,  2,3, \cdots. \label{decay-1d-b}
\end{align}\end{subequations}
This is estimate (2.13) in \cite{LZ}, whose proof follows from the ODE analyses.

{\em Step 1}. In this step, we prove \ef{1.31} for $j=0$, that is,
\begin{align}\label{1.31-1}
\mathscr{E}_0(t)+\int_0^t \mathscr{D}_0(s)ds \le C(\ga, M) \mathscr{E}_0(0).
\end{align}

Multiply \ef{eluerpert} by $w_t$ and integrate the resulting equation over $\mathcal{I}$ to get
 \begin{align}\label{12.31-1}
\f{d}{dt} \mathfrak{E}_0(t)
 +
 \int_\mathcal{I}  \bar\rho_0w_t^2 dx =- \tilde{\eta}_{xt}  \int_\mathcal{I}  \bar\rho_0^\ga R dx ,
\end{align}
where $R$ is defined by \ef{date2.10}, and
 \begin{align*}
\mathfrak{E}_0(t) =\int_\mathcal{I} \lt\{ \frac{1}{2}\bar\rho_0w_{t}^2   + \frac{1}{\ga-1} \bar\rho_0^\ga \lt[(\tilde{\eta}_x+w_x)^{1-\ga } -\tilde{\eta}_x^{1-\ga }-(1-\ga)\tilde{\eta}_x^{-\ga}   w_{x}\rt] \rt\}dx.
\end{align*}
This, together with the Taylor expansion, \ef{decay-1d-a} and \ef{20lem31est}, implies that
\begin{align}\label{12.31-2}
(C(\ga))^{-1}  { {E}}_0(t) \le  \mathfrak{E}_0(t) dx\le C(\ga) { {E}}_0(t) ,\ \
 R \ge C(\ga)^{-1} \tilde{\eta}_x^{-\ga-2}    w_x^2 \ge 0 ,
\end{align}
where
$$E_0(t)=\int_\mathcal{I}  \bar\rho_0    |w_t|^2    dx + (1+ t)^{-1} \int_\mathcal{I}     \bar\rho_0^\ga  | w_x|^2   dx  .$$
Integrate \ef{12.31-1} over $[0,t]$, and use \ef{decay-1d-a} and \ef{12.31-2} to achieve
\be\label{12.31-3}
 { {E}}_0(t)+\int_0^t \int_\mathcal{I}  \bar\rho_0w_s^2 dx ds \le C(\ga) { {E}}_0(0).
\ee

Multiply \ef{12.30} by ${\tilde\eta_x}^{\ga+1} w_t$ and integrate the resulting equation over $\mathcal{I}$ to give
\begin{align}\label{1.2-1}
\f{d}{dt}\int_\mathcal{I} \lt( \frac{1}{2}{\tilde\eta_x}^{\ga+1}  \bar\rho_0 w_t^2 + \bar\rho_0^\ga\lt(\frac{\ga}{2} w_x^2 +  {G}\rt)\rt)   dx
 +
 {\tilde\eta_x}^{\ga+1} \int_\mathcal{I}  \bar\rho_0w_t^2 dx \notag  \\
 =- \tilde{\eta}_{xt}  \int_\mathcal{I}  \bar\rho_0^\ga   {H} dx +  \frac{1}{2} \lt( {\tilde\eta_x}^{\ga+1} \rt)_t \int_\mathcal{I}   \bar\rho_0 w_t^2   dx,
\end{align}
where
\begin{align*}
& {G}(t)=\frac{1}{\ga-1} {\tilde\eta_x}^{\ga+1} (\tilde\eta_x+w_x)^{1-\ga}
- \frac{1}{\ga-1}  {\tilde\eta_x}^2 +  {\tilde\eta_x} w_x - \frac{\ga}{2} w_x^2, \\
& {H}(t)=\frac{\ga+1}{\ga-1}{\tilde\eta_x}^{\ga} (\tilde\eta_x+w_x)^{1-\ga} -  {\tilde\eta_x}^{\ga+1} (\tilde\eta_x+w_x)^{ -\ga} -\frac{2}{\ga-1} {\tilde\eta_x} + w_x.
 \end{align*}
It follows from the Taylor expansion,  \ef{decay-1d-a} and \ef{20lem31est} that
\begin{align}\label{1.3-1}
| {G}(t)|\le C(\ga) (1+t)^{-\frac{1}{\ga+1}} |w_x|^3
 \ \ {\rm and} \ \  | {H}(t)| \le C(\ga) (1+t)^{-\frac{2}{\ga+1}} |w_x|^3,
 \end{align}
then we integrate \ef{1.2-1} over $[0,t]$, and use \ef{decay-1d} and \ef{20lem31est} to obtain
\begin{align*}
& (1+t) { {E}}_0(t)+\int_0^t (1+s) \int_\mathcal{I}  \bar\rho_0w_s^2 dx ds \le  C(\ga )  { {E}}_0(0)   \\
& + C(\ga,M)\epsilon_0 \int_0^t (1+s)^{-1-\frac{1}{\ga+1}} \int_\mathcal{I} \bar\rho_0^\ga w_x^2 dx ds
 + C(\ga)\int_0^t   \int_\mathcal{I}  \bar\rho_0w_s^2 dx ds.
\end{align*}
This, together with  \ef{12.31-3},  means that
\begin{align}\label{1.3-2}
 (1+t) { {E}}_0(t)+\int_0^t (1+s) \int_\mathcal{I}  \bar\rho_0w_s^2 dx ds
 \le  C(\ga)  { {E}}_0(0).
\end{align}

Multiply  \ef{eluerpert} by $w$ and integrate the resulting equation over $\mathcal{I}$ to get
\begin{align}\label{1.3-3}
 \f{1}{2}\frac{d}{dt} \int_\mathcal{I} \bar\rho_0\lt( w^2 + 2ww_t \rt)  dx  -
 \int_\mathcal{I}  \bar\rho_0^\ga\lt[ (\tilde{\eta}_{x}+w_x)^{-\ga} -\tilde{\eta}_{x}^{-\ga}   \rt] w_{x}   dx =     \int_\mathcal{I}  \bar\rho_0w_t^2 dx .
\end{align}
It follows from the Taylor expansion, \ef{decay-1d-a} and \ef{20lem31est} that
 $$-\lt[ (\tilde{\eta}_{x}+w_x)^{-\ga} -\tilde{\eta}_{x}^{-\ga}   \rt] w_{x}
 \ge C(\ga)^{-1} (1+t)^{-1}  w_{x}^2,$$
 then we integrate \ef{1.3-3} over $[0,t]$, and use the Cauchy inequality and \ef{12.31-3} to obtain
\begin{align}
  \int_\mathcal{I} \lt(\bar\rho_0 w^2 \rt)(x,t)  dx + \int_0^t (1+s)^{-1}  \int_\mathcal{I}  \bar\rho_0^\ga  w_{x}^2   dx ds  \le C(\ga) \mathscr{E}_0(0). \label{1.4-1}
\end{align}

Finally, we multiply \ef{1.2-1} by $\ln(1+t)$ and integrate the resulting equation over $[0,t]$ to achieve
\begin{align}
&(1+t)\ln(1+t) { {E}}_0(t)+\int_0^t (1+s) \ln(1+s) \int_\mathcal{I}  \bar\rho_0w_s^2 dx ds \notag \\
 \le & C(\ga)  { {E}}_0(0)
 +C(\ga) \int_0^t (1+s)^{-1} \int_\mathcal{I}  \bar\rho_0^\ga w_x^2 dx ds \notag \\
& + C(\ga,M)\epsilon_0 \int_0^t (1+s)^{-1-\frac{1}{\ga+1}} \lt(1+ \ln(1+s) \rt) \int_\mathcal{I}  \bar\rho_0^\ga w_x^2 dx ds\notag \\
& + C(\ga)\int_0^t \lt( 1+ \ln(1+s)  \rt)\int_\mathcal{I}  \bar\rho_0w_s^2 dx ds \notag\\
\le & C(\ga) {{E}}_0(0)  +C(\ga,M) \int_0^t (1+s)^{-1} \int_\mathcal{I}  \bar\rho_0^\ga w_x^2 dx ds \notag \\
& + C(\ga)\int_0^t (1+s)\int_\mathcal{I}  \bar\rho_0w_s^2 dx ds \le C(\ga,M){\mathscr{E}}_0(0), \label{d3.11-1}
\end{align}
where the first inequality follows from \ef{1.3-1},  \ef{decay-1d} and \ef{20lem31est}; the second from $\ln(1+t)\le 1+t$ and $\ln(1+t)\le C(\ga)(1+t)^{1/(\ga+1)}$; and the last from \ef{1.3-2} and \ef{1.4-1}.  This, together with  \ef{1.4-1}, proves \ef{1.31-1}.

{\em Step 2}.  In this step, we use the mathematical induction to prove \ef{1.31} for $j \ge 1$. It follows from \ef{1.31-1} that \ef{1.31} holds for $j=0$. For $1\le i \le 4+ [\alpha]$, we make the hypothesis that \ef{1.31} holds for $j=0, 1, \cdots, i-1$, that is,
\begin{align}\label{1.31-2}
\mathscr{E}_j(t)+\int_0^t \mathscr{D}_j(s)ds \le C(\ga, M)\sum_{0\le m\le j}\mathscr{E}_m(0), \ \  j=0,1,\cdots,i-1.
\end{align}
It is enough to prove \ef{1.31} holds for $j=i$ under the hypothesis \ef{1.31-2}.

Multiply \ef{12.30} by $\bar\rho_0^{-1}$, and
take $\pl_x^i$  onto  the resulting equation to get
 \begin{align}\label{1.14}
\pl_x^i w_{tt}   +
  \pl_x^i w_t
  -
   \varsigma^{-\alpha-i} \lt( \varsigma^{1+\alpha+i} \pl_x^{i} \lt(  \ga {\tilde\eta_x}^{-\ga-1} w_x  -{R} \rt)\rt)_x \notag\\
   +m_i \pl_x^{i-1}\lt(  \ga {\tilde\eta_x}^{-\ga-1} w_x -{R}  \rt) =0, \ \  i\ge 1,
   \end{align}
where $\varsigma=\bar\rho_0^{\ga-1}=A_1-B_1x^2$  and $m_i=  i(1+2\alpha +i)B_1 >0$. Note that for $i\ge 1$,
\begin{align}\label{date2.4}
\pl_x^i {R} =    \ga (\tilde\eta_x^{-\ga-1}-(\tilde\eta_x + w_x)^{-\ga-1} )  \pl_x^{i } w_x - {N}_i,
\end{align}
where
\begin{align}\label{20.2.17}
{N}_i=   \ga   \lt( \pl_x^{i-1}( (\tilde\eta_x+ w_x)^{-\ga-1} \pl_x  w_x ) -(\tilde\eta_x+  w_x)^{-\ga-1} \pl_x^{i } w_x  \rt).
\end{align}
Then, we  multiply \ef{1.14} by $\varsigma^{\alpha+i}\pl_x^i w_t$,  integrate the resulting equation over $\mathcal{I}$, and use \ef{decay-1d-a} which implies  $({\tilde\eta_x}^{-\ga-1})_t\le 0$ to obtain
\begin{align}
& \frac{d}{dt}E_i(t)
 +\int_\mathcal{I} \varsigma^{\alpha+i}|\pl_x^i w_t|^2dx \notag\\
=&\frac{\ga}{2} ({\tilde\eta_x}^{-\ga-1})_t
\int_\mathcal{I} \varsigma^{\alpha+i} \lt( \varsigma|\pl_x^{i } w_x |^2 + m_i |\pl_x^i w|^2 \rt)dx + Z_i \le Z_i, \label{date21}
\end{align}
where
\begin{align*}
&E_i(t)=
\frac{1}{2} \int_\mathcal{I} \varsigma^{\alpha+i} \lt( |\pl_x^i w_t|^2 + \ga {\tilde\eta_x}^{-\ga-1} \lt( \varsigma
|\pl_x^{i } w_x |^2  + m_i |\pl_x^{i} w|^2 \rt) \rt)dx+ Q_i(t),\\
 & Q_i(t)=   \int_\mathcal{I} \varsigma^{1+\alpha+i}  \lt(   \frac{\ga}{2}((\tilde\eta_x + w_x)^{-\ga-1}- \tilde\eta_x^{-\ga-1}  )  |\pl_x^{i} w_x |^2  +   {N}_i  \pl_x^{i } w_x  \rt) dx,\\
 & Z_i(t)=   \frac{\ga}{2}  \int_\mathcal{I} \varsigma^{1+\alpha+i}  \pl_t\lt(   (\tilde\eta_x + w_x)^{-\ga-1}- \tilde\eta_x^{-\ga-1}  \rt)  |\pl_x^{i} w_x|^2  dx \\
&  \ \ +  \int_\mathcal{I} \varsigma^{1+\alpha+i} (\pl_t   {N}_i) \pl_x^{i } w_x  dx+m_i \int_\mathcal{I}  \varsigma^{\alpha+i}  ( \pl_x^{i-1} {R} ) \pl_x^i w_t dx.
\end{align*}
Similarly, we multiply \ef{1.14} by $\tilde\eta_x^{\ga+1}\varsigma^{\alpha+i}\pl_x^i w_t$,   integrate the resulting equation over $\mathcal{I}$, and use \ef{decay-1d}  to give
\begin{align}\label{2-1}
& \frac{d}{dt}\lt(\tilde\eta_x^{\ga+1} E_i(t)\rt)+  \tilde\eta_x^{\ga+1}\int_\mathcal{I} \varsigma^{\alpha+i}|\pl_x^i w_t|^2dx= Y_i \notag\\
& +\frac{1}{2} \lt(\tilde\eta_x^{\ga+1} \rt)_t \int_\mathcal{I} \varsigma^{\alpha+i}|\pl_x^i w_t|^2dx \le Y_i + c_1(\ga) \int_\mathcal{I} \varsigma^{\alpha+i}|\pl_x^i w_t|^2dx
\end{align}
for a certain positive constant $c_1$ depending only on $\ga$,
where
\begin{align*}
& Y_i(t)=     \frac{\ga}{2}  \int_\mathcal{I} \varsigma^{1+\alpha+i}  \pl_t \lt(   (1+ \tilde\eta_x^{-1}   w_x)^{-\ga-1}  \rt) |\pl_x^{i} w_x|^2  dx \\
& +  \int_\mathcal{I} \varsigma^{1+\alpha+i} (\pl_t(   \tilde\eta_x^{\ga+1}{N}_i)) \pl_x^{i } w_x  dx +
  m_i\tilde\eta_x^{\ga+1} \int_\mathcal{I}  \varsigma^{\alpha+i}  ( \pl_x^{i-1} {R} ) \pl_x^i w_t dx.
\end{align*}
Multiply \ef{2-1} by $\ln(1+t)$, and use \ef{decay-1d-a} to achieve
\begin{align}\label{2-1-1}
&\frac{d}{dt}\lt(\ln(1+t)\tilde\eta_x^{\ga+1}
 E_i(t)\rt)+  \ln(1+t)\tilde\eta_x^{\ga+1}\int_\mathcal{I} \varsigma^{\alpha+i}|\pl_x^i w_t|^2dx\notag\\
\le & c_2(\ga)
E_i + \ln(1+t)Y_i + c_1(\ga)   \tilde\eta_x^{\ga+1}  \int_\mathcal{I} \varsigma^{\alpha+i}|\pl_x^i w_t|^2dx
\end{align}
for a certain positive constant $c_2$ depending only on $\ga$,
where $(1+t)^{-1}\tilde\eta_x^{\ga+1}(t)\le c_2(\ga)$ and $\ln(1+t)\le 1+t\le \tilde\eta_x^{\ga+1}(t)$ which follow from \ef{decay-1d-a}, have been used to derive \ef{2-1-1}.
To control the first term on the right hand side of \ef{2-1-1}, we multiply \ef{1.14} by $ \varsigma^{\alpha+i}\pl_x^i w $, and integrate the resulting equation over $\mathcal{I}$ to obtain
\begin{align}\label{2-1-2}
 \frac{d}{dt}\int_\mathcal{I} \varsigma^{\alpha+i}\lt(\frac{1}{2}|\pl_x^i w|^2 dx+   (\pl_x^i w)\pl_x^i w_t\rt)dx
 +  \ga {\tilde\eta_x}^{-\ga-1} \int_\mathcal{I}\varsigma^{\alpha+i}   ( \varsigma
|\pl_x^{i } w_x |^2 \notag\\
   + m_i |\pl_x^{i} w|^2  ) dx
=   \int_\mathcal{I} \varsigma^{\alpha+i}|\pl_x^i w_t|^2dx +   X_i,
\end{align}
where
$$X_i(t)= \int_\mathcal{I} \varsigma^{\alpha+i}   \lt( \varsigma
( \pl_x^{i } {R} ) \pl_x^{i } w_x    + m_i ( \pl_x^{i-1} {R} )  \pl_x^{i} w  \rt) dx . $$
Finally, it produces from the summation of
$(1+c_1(\ga))(\ef{2-1}
+c_1(\ga)\ef{date21})$, $\ef{2-1-1}$ and $(1+2^{-1}c_2(\ga))
(\ef{2-1-2}+2\ef{date21})$ that
\begin{align}\label{2-5-1}
\frac{d}{dt}\mathfrak{E}_i(t)+\mathfrak{D}_i(t)
\le C(\ga) W_i(t),
\end{align}
where
\begin{align*}
 \mathfrak{E}_i(t)= & \lt(1+\frac{1}{2}c_2(\ga)\rt)\lt(\int_\mathcal{I} \varsigma^{\alpha+i}\lt(\frac{1}{2}|\pl_x^i w|^2 dx+   (\pl_x^i w)\pl_x^i w_t\rt)dx+2E_i(t)\rt)\\
&   +\lt((1+c_1(\ga))(\tilde\eta_x^{\ga+1}
+c_1(\ga))+\ln(1+t)
\tilde\eta_x^{\ga+1}\rt)E_i(t),\\
 \mathfrak{D}_i(t)= & \lt(1+(1+\ln(1+t))
\tilde\eta_x^{\ga+1}\rt)\int_\mathcal{I} \varsigma^{\alpha+i}|\pl_x^i w_t|^2dx \\
& + \ga {\tilde\eta_x}^{-\ga-1} \int_\mathcal{I}\varsigma^{\alpha+i}   ( \varsigma
|\pl_x^{i } w_x |^2
   + m_i |\pl_x^{i} w|^2  ) dx , \\
   W_i(t)=& |Q_i|+|X_i|+(1+\ln(1+t))|Y_i|
+|Z_i| .
   \end{align*}

Next, we analyze \ef{2-5-1}. Due to the Taylor expansion, the Cauchy inequality, \ef{decay-1d}, \ef{20lem31est}, \ef{date2.10} and \ef{date2.4}, we have that for any constant $\epsilon>0$,
\begin{subequations}\label{2-5-3}\begin{align}
&|Q_i(t)|\le   \lt( \epsilon+ C(\ga,M)\epsilon_0 \rt) (1+t)^{-1}
\int_\mathcal{I}\varsigma^{1+\alpha+i}|
\pl_x^{i}w_x|^2dx  +\epsilon^{-1}H_i(t), \label{2-5-3-a}\\
& W_i(t)\le \lt(\epsilon +C(\ga,M) \epsilon_0   \rt)\mathfrak{D}_i(t)
+\epsilon^{-1}\lt(C(\ga) G_{i }
+C(\ga,M) \epsilon_0^2\mathscr{D}_{i-1}  \rt)(t),\label{2-5-3-b}
\end{align}\end{subequations}

where
\begin{subequations}\label{20-1}\begin{align}
 H_i(t)=&(1+t)\int_\mathcal{I}\varsigma^{1+\alpha+i}|
N_i|^2dx,\label{20.2.19-3}\\
 G_{i }(t)= & (\ln(1+t)+1)^2H_i(t)  + (\ln(1+t)+1) H_{i-1}(t)
 \notag \\ &  +(1+t)^3(\ln(1+t)+1)^2\int_\mathcal{I}\varsigma^{1+\alpha+i}|
\pl_t{N}_i|^2dx\label{20.2.19-4}
\end{align}\end{subequations}
with $N_0(t)=0$. This, together with the Cauchy inequality and \ef{decay-1d-a}, implies that \begin{subequations}\label{2-5-2}\begin{align}
& \mathfrak{E}_i(t)\ge C(\ga)^{-1}\mathscr{E}_i(t)
 -C(\ga)(1+t)(\ln(1+t)+1)H_i(t),
  \\
&    \mathfrak{E}_i(t)\le C(\ga) ( \mathscr{E}_i +  \mathscr{E}_{i-1})(t)
 +C(\ga)(1+t)(\ln(1+t)+1) H_i(t),\\
&  \mathfrak{D}_i(t) \ge C(\ga)^{-1}\mathscr{D}_i(t).
\end{align}\end{subequations}

To estimate $H_i$ and $G_i$, we
note that $N_i$ defined in \ef{20.2.17} can be rewritten as
\begin{align*}
{N}_i
= \ga \tilde\eta_x^{-\ga-1}\sum_{1\le m\le i-1} C_{i-1}^m  (\pl_x^m (1+\tilde\eta_x^{-1} w_x)^{-\ga-1} ) \pl_x^{i-m}  w_x,
\end{align*}
where $C_{i-1}^m=\frac{m!}{(i-1)!(i-1-m)!}$
and $\sum_{1\le m\le i-1}$ should be understood as zero when $i=0,1$; and
 \begin{align*}
 |\pl_x^k (1+ \tilde\eta_x^{-1}w_x)^{-\ga-1}|+ |\pl_x^k (1+ \tilde\eta_x^{-1}w_x)^{-\ga-2}|
 \le C(\ga,k) J_k,  \ \ k\ge 0,
\end{align*}
due to \ef{decay-1d-a} and \ef{20lem31est}, where  $J_k$ $(k\ge 0)$ are defined  inductively by
\begin{align*}
J_0=1 ,   \ \
J_k=\tilde\eta_x^{-1} \lt(|\pl_x^k w_x|+\sum_{1\le m\le k-1}J_m |\pl_x^{k-m} w_x|\rt).
\end{align*}
Then,  we have
$N_i=0$ for $i=0,1$, and for $2\le i \le 4+[\alpha]$,
\begin{subequations}\label{20-2}\begin{align}
&|N_i|\le C(\ga)(1+t)^{-1}\sum_{1\le m\le i-1} J_m |\pl_x^{i-m}w_x|, \label{20.2.19-1}\\
&|\pl_t N_i| \le  C(\ga) (1+t)^{-2}
\sum_{1\le m\le i-1} J_m |\pl_x^{i-m}w_x| + C(\ga) (1+t)^{-1-\frac{1}{\ga+1}}\notag\\
& \qquad \times
\sum_{1\le m\le i-1,\ 0\le \iota\le m}J_\iota\lt((1+t)^{-1}\pl_x^{m-\iota}w_x+\pl_x^{m+1-\iota}w_t\rt)  |\pl_x^{i-m}w_x|\notag\\
& \qquad + C(\ga)(1+t)^{-1} \sum_{1\le m\le i-1}J_m |\pl_x^{i+1-m}w_t|, \label{20.2.19-2}
\end{align}\end{subequations}
with the aid of \ef{decay-1d}.
It follows from \ef{decay-1d-a} and \ef{20lem31est} that for $1\le m\le i-1$,
\begin{align*}
& (1+t)^{ \frac{2}{\ga+1}} \varsigma^{ {2m-1} }J_m^2 \\
\le &
 C(\ga)
\lt(\varsigma^{ {2m-1} }|\pl_x^m w_x|^2+\sum_{1\le \iota\le m-1}\varsigma^{  {2m-1}  }J_\iota^2 |\pl_x^{m-\iota} w_x|^2\rt) \\
\le & C(\ga,M)\ea_0^2 + C(\ga)\sum_{1\le \iota\le m-1}\varsigma^{  {2\iota-1}  }J_\iota^2  \varsigma^{  {2(m-\iota)}  }  |\pl_x^{m-\iota} w_x|^2\\
\le & C(\ga,M)\ea_0^2 + C(\ga,M) \ea_0^2 (1+t)^{- \frac{2}{\ga+1}}  \sum_{1\le \iota\le m-1} (1+t)^{ \frac{2}{\ga+1}} \varsigma^{  {2\iota-1}  }J_\iota^2\\
\le  & C(\ga,M)\ea_0^2,
\end{align*}
which, together with \ef{hardy2019-1d}, \ef{20.2.19-3} and \ef{20.2.19-1}, implies that for any positive integer $k$,
\begin{align}
& (\ln(1+t)+1)^k H_i(t) \notag\\
 \le & C(\ga)(\ln(1+t)+1)^k (1+t)^{-1}\sum_{1\le m\le i-1}\int_\mathcal{I}\varsigma^{1+\alpha+i} J_m^2  |\pl_x^{i-m}w_x|^2dx \notag\\
 \le &
 C(\ga,M,k)\ea_0^2 (1+t)^{-1} \sum_{1\le m\le i-1} \int_\mathcal{I} \varsigma^{2+\alpha+i-2m} |\pl_x^{i-m}w_x|^2dx\notag \\
\le & C(\ga,M,k)\ea_0^2 (1+t)^{-1} \sum_{1\le m\le i-1} \int_\mathcal{I} \varsigma^{ \alpha+m+1} |\pl_x^{m}w_x|^2dx. \label{20March22}
\end{align}
Indeed, the following estimate has been used to derive the last inequality of
\ef{20March22}.
\begin{align*}
&
 \sum_{1\le m\le i-1} \int_\mathcal{I} \varsigma^{2+\alpha+i-2m} |\pl_x^{i-m}w_x|^2dx
\\
\le &  \sum_{1\le m\le i-1} \int_\mathcal{I} \varsigma^{2+\alpha+i-2m+2} (|\pl_x^{i-m}w_x|^2+|\pl_x^{i-m+1}w_x|^2)dx\\
\le &  \sum_{1\le m\le i-1} \int_\mathcal{I} \varsigma^{2+\alpha+i-2m+2(m-1)}\sum_{0\le \iota\le m-1} |\pl_x^{i-m+\iota}w_x|^2 dx\\
= &  \sum_{1\le m\le i-1} \int_\mathcal{I} \varsigma^{ \alpha+i}\sum_{i-m\le \iota\le i-1} |\pl_x^{\iota}w_x|^2 dx,
\end{align*}
where the first inequality follows from  \ef{hardy2019-1d} and $2+\alpha+i-2m\ge 4+\alpha-i\ge 0$, and the second from the repeat  use of \ef{hardy2019-1d}.
It follows from \ef{20March22} that
\begin{align}
 (1+t)(\ln(1+t)+1)H_i(t)
\le   C(\ga,M)\ea_0^2 \sum_{1\le m\le i-1} \mathscr{E}_m(t). \label{20.2.19-a}
\end{align}
Similarly, we can use \ef{20-1} and \ef{20-2} to get
\begin{align}
G_i(t)\le C(\ga, M)\epsilon_0^2
 \sum_{0\le m\le i  } \mathscr{D}_m(t). \label{20-2-17}
\end{align}
For example, we examine the most difficult term. It holds that for any positive integer $k$,
\begin{align}
& (1+t)^3(\ln(1+t)+1)^k \int_\mathcal{I}
\varsigma^{1+\alpha+i} \lt((1+t)^{-1} \sum_{1\le m\le i-1}J_m |\pl_x^{i+1-m}w_t|\rt)^2 dx\notag \\
&\le  C(\ga) (1+t)(\ln(1+t)+1)^k \sum_{1\le m\le i-1}\int_\mathcal{I}\varsigma^{1+\alpha+i} J_m^2  |\pl_x^{i+1-m}w_t|^2dx
\notag  \\
 &\le    C(\ga,M,k)\ea_0^2  (1+t) \sum_{1\le m\le i-1} \int_\mathcal{I} \varsigma^{2+\alpha+i-2m}
 |\pl_x^{i+1-m}w_t|^2dx\notag \\
&\le   C(\ga,M,k)\ea_0^2 (1+t) \sum_{1\le m\le i-1} \int_\mathcal{I} \varsigma^{ \alpha+i }\sum_{i+1-m\le \iota\le i} |\pl_x^{ \iota}w_t|^2 dx\notag \\
& \le   C(\ga, M,k)\epsilon_0^2
 \sum_{2\le m\le i  } \mathscr{D}_m(t). \label{20March22-1}
\end{align}

Finally, we integrate \ef{2-5-1} over $[0,t]$, and use \ef{1.31-2}, \ef{2-5-3-b}, \ef{2-5-2}, \ef{20.2.19-a} and \ef{20-2-17} to prove that \ef{1.31-2} holds for $j=i$.
This finishes the proof of Lemma  \ref{lem2020}.

\subsubsection{Proof of Lemma \ref{lem-d3.8}}
We use the mathematical induction to prove  \ef{new12-4} for $j=0$. The key is to show the zeroth order estimate which determines the time decay rate.

{\em Step 1}.
In this step, we prove that \ef{new12-4} holds for $j=0$, that is,
 \begin{align}
&\lt(\lt( \ln(1+t) \rt)^q +1 \rt)\mathcal{E}_0(t; w-\vartheta_0)  \notag \\
&\quad +  \int_0^t \lt( \lt( \ln(1+s) \rt)^q +1 \rt) \mathcal{D}_{0}(s) ds
  \le C(\ga, M,q)   \mathcal{E} (0).\label{d3.23-1}
\end{align}

Let $\vartheta(t)$ be the center of mass, then it follows from \ef{Feb25-3} that
$$\vartheta(t)=\frac{1}{M}\int_{\Oa(t)} \eta \rho(\eta,t)d\eta=\frac{1}{M}\int_\mathcal{I}  \bar\rho_0(x) \eta(x,t)dx,$$
We integrate  \ef{equa} over $\mathcal{I}$ to obtain  $\vartheta''(t) +\vartheta'(t) =0$, which implies
\begin{align}\label{Feb28-2}
& \vartheta(t)=   \vartheta(0)+(1-e^{-t})\vartheta'(0)
=\vartheta_0-   e^{-t} \frac{1}{M}  \int_{\Oa(0)} \rho_0(x) u_0(x) dx,
\end{align}
where $\vartheta_0$ is defined by \ef{d3.12-1}.
It follows from $\int_{\mathcal{I}} x\bar\rho_0(x) dx =0 $ that
\begin{align*}
 \int_\mathcal{I} \bar\rho_0(x) \lt(w(x,t)-   \vartheta(t) \rt) dx
=
\int_\mathcal{I} \bar\rho_0(x)  \eta(x,t) dx-    M \vartheta(t)  =0.
\end{align*}
Then, there exists $x=b(t)\in \mathcal{I}$ such that
 $\bar\rho_0(b(t)) \lt(w(b(t),t)-   \vartheta(t) \rt)=0$, which means
$w(b(t),t)=\vartheta(t)$ and then
$$
 w(x,t)-  \vartheta(t) = \int_{b(t)}^x \pl_y \lt( w(y,t)-  \vartheta(t) \rt) dy  = \int_{b(t)}^x w_y(y,t)dy.
$$
This, together with \ef{Feb28-2}, gives
\begin{align} \label{d3.10}
 \lt| w(x,t)- \vartheta_0 \rt|
 \le  C(\ga, M)e^{-t} \sqrt{\mathcal{E}(0)} + \lt|  \int_{b(t)}^x   w_y(y,t)dy \rt|
\end{align}
for $b(t)\in \mathcal{I}$,
which implies
\begin{align}
&\int_{\mathcal{I}} (w(x,t)-\vartheta_0)^2 dx \le C(\ga, M)  e^{-2t} \mathcal{E}(0)  + C(\ga, M)  \int_{\mathcal{I}} w_x^2 dx . \label{d3.12-2}
\end{align}
Indeed, the following estimate has been used to derive \ef{d3.10}.
\begin{align*}
  \lt|\frac{1}{M} \int_{\Oa(0)}  \rho_0  u_0  dx \rt| =\lt|\frac{1}{M} \int_{ \mathcal{I}} \bar \rho_0  \eta_t\big|_{t=0} dx \rt|
= \lt| \frac{1}{M}\int_{ \mathcal{I}} \bar \rho_0  w_t\big|_{t=0} dx \rt|  \\
  \le   \lt\|w_t(\cdot, 0) \rt\|_{L^\iy(\mathcal{I})} \le C(\ga, M) \sqrt{\mathcal{E}(0)},
\end{align*}
where the second equality follows from $\int_{\mathcal{I}} x\bar\rho_0(x) dx =0 $,  and the last inequality from  \ef{20lem31est}.
It follows from \ef{wsv}
and   $\bar\rho_0^{\ga-1}=A-Bx^2$ that
\begin{align}
&\|w_x\|_{L^2(\mathcal{I})}^2\le \|w_x  \|_{H^{ (1+[\alpha]-\alpha)/{2}}(\mathcal{I})}^2
\le   C(\ga,M)\|w_x  \|_{H^{\alpha+3+[\alpha],\ {2+[\alpha]}}(\mathcal{I})}^2 \notag\\
 & =  C(\ga,M) \sum_{0\le i \le 2+[\alpha]} \int_\mathcal{I} \bar\rho_0^{(\ga-1)(\alpha+3+[\alpha])} |\pl_x^i w_x|^2 dx
 \notag\\
& \le   C(\ga,M) \sum_{0 \le i \le 2+[\alpha]} \int_\mathcal{I}  \bar\rho_0^{\ga+i(\ga-1)} |\pl_x^i w_x |^2 dx ,\notag
\end{align}
which, together with
\ef{d3.12-2}, means
\begin{align}
&\int_{\mathcal{I}} (w-\vartheta_0)^2  dx   \le
 C(\ga, M)  e^{-2t} \mathcal{E}(0)\notag \\
&\quad +C(\ga, M)   \sum_{0 \le i \le 2+[\alpha]}  \int_\mathcal{I}  \bar\rho_0^{\ga+i(\ga-1)} |\pl_x^i w_x |^2 dx.\label{d3.10-2}
\end{align}

Let $k$ be any given positive integer and $\varsigma=\bar\rho_0^{\ga-1}$. We integrate the product of  $\lt(\ln(1+t)\rt)^k $ and \ef{2-1-2} over $[0,t]$, and use the Cauchy inequality, \ef{decay-1d-a} and \ef{1.31} to obtain that
\begin{align}\label{3.22-1}
&\lt(\ln(1+t)\rt)^k \int_\mathcal{I} \varsigma^{\alpha+i} |\pl_x^i w|^2 dx \notag\\
&+\int_0^t (1+s)^{-1} \lt(\ln(1+s)\rt)^k
  \int_\mathcal{I}\varsigma^{\alpha+i}   ( \varsigma
|\pl_x^{i } w_x |^2
   +   |\pl_x^{i} w|^2  ) dx ds\notag\\
&   \le C(\ga,M,k) \sum_{0\le m\le i} \mathscr{E}_m(0)
 +C(\ga)  \int_0^t \lt(\ln(1+s)\rt)^k   X_i(s) ds\notag\\
& + C(\ga) k \int_0^t \lt(\ln(1+s)\rt)^{k-1} (1+s)^{-1}
  \int_\mathcal{I}\varsigma^{\alpha+i}    |\pl_x^{i} w|^2   dx ds
\end{align}
for $1\le i \le 4+[\alpha]$.
It follows from \ef{date2.4}, \ef{decay-1d-a}, \ef{20lem31est}, \ef{20.2.19-3} and \ef{20March22} that
\begin{align}\label{3.22-4}
&\lt(\ln(1+t)\rt)^k   X_i(t)
\le \lt(C(\ga,M,k)\ea_0+1\rt)(1+t)^{-1}\int_\mathcal{I}\varsigma^{\alpha+i}   ( \varsigma
|\pl_x^{i } w_x |^2
 \notag \\
&   +   |\pl_x^{i} w|^2  ) dx
   +\lt(\ln(1+t)\rt)^k \lt(H_i(t)+H_{i-1}(t)\rt)
   \le C(\ga, M, k) \sum_{0\le m\le i} \mathscr{D}_m(t),
\end{align}
which, together with \ef{1.31}, \ef{3.22-1} and the mathematical induction, gives that
\begin{align}\label{3.22-2}
& \int_0^t (1+s)^{-1} \lt(\ln(1+s)\rt)^q
  \int_\mathcal{I}\varsigma^{\alpha+i}   ( \varsigma
|\pl_x^{i } w_x |^2
   +   |\pl_x^{i} w|^2  ) dx ds   \notag\\
& + \lt(\ln(1+t)\rt)^q \int_\mathcal{I} \varsigma^{\alpha+i} |\pl_x^i w|^2 dx \le C(\ga,M,q) \sum_{0\le m\le i} \mathscr{E}_m(0)
\end{align}
for $1\le i \le 4+[\alpha]$.
Based on \ef{3.22-2} and \ef{1.31}, we
integrate the product of  $\lt(\ln(1+t)\rt)^q $ and \ef{2-1} over $[0,t]$ to get that for $1\le i \le 4+[\alpha]$,
\begin{align*}
& (\ln(1+t))^q
\int_\mathcal{I} \varsigma^{\alpha+i} \lt( (1+t)
|\pl_x^i w_t|^2 +  \varsigma
|\pl_x^{i } w_x |^2  +  |\pl_x^{i} w|^2  \rt)
dx
\notag
\\
&+\int_0^t (1+s) (\ln(1+s))^q   \int_\mathcal{I} \varsigma^{\alpha+i} |\pl_x^i w_s|^2 dx ds
 \le  C(\ga,M,q) \sum_{0\le m\le i} \mathscr{E}_m(0),
\end{align*}
where the derivation follows similarly from those of \ef{20March22-1} and \ef{3.22-4}. This, together with
\ef{d3.10-2} and \ef{1.31}, implies that
\begin{align}
&\lt(\lt( \ln(1+t) \rt)^q +1 \rt)\int_{\mathcal{I}} (w-\vartheta_0)^2  dx   \le
 C(\ga, M,q)    \mathcal{E}(0).\label{d3.23}
\end{align}

In a similar way to deriving \ef{d3.11-1},
we integrate the product of  $\lt(\ln(1+t)\rt)^q $ and \ef{1.2-1} over $[0,t]$, and use \ef{3.22-2} to achieve
\begin{align*}
& (\ln(1+t))^q \int_\mathcal{I} \lt( (1+t) \bar\rho_0    |w_t|^2    +  \bar\rho_0^\ga  | w_x|^2  \rt) dx \notag
\\
&+\int_0^t (1+s) (\ln(1+s))^q \int_\mathcal{I}  \bar\rho_0w_s^2 dx ds
 \le  C(\ga,M,q)  \sum_{0\le m\le 1} \mathscr{E}_m(0),
\end{align*}
which, together with \ef{1.31}, \ef{d3.23}, \ef{3.22-2}, proves \ef{d3.23-1}.

{\em Step 2}.
In this step, we use the mathematical induction to prove \ef{new12-4} for $j \ge 1$. It follows from \ef{d3.23-1} that \ef{new12-4} holds for $j=0$. For $1\le k\le 4+ [\alpha]$, we make the hypothesis that \ef{new12-4} holds for $j=0, 1, \cdots, k-1$, that is,
 \begin{align}
\lt(\lt( \ln(1+t) \rt)^q +1 \rt)\mathcal{E}_j(t;w-\vartheta_0) +  \int_0^t \lt( \lt( \ln(1+s) \rt)^q +1 \rt) \mathcal{D}_{j}(s) ds  \notag \\
  \le C(\ga, M,q )   \mathcal{E} (0), \ \  \ \  j=0 ,1 ,2, \cdots,  k-1.  \label{12-4-j-20}
\end{align}
It is enough to prove \ef{new12-4} holds for $j=k$ under the hypothesis \ef{12-4-j-20}. For simplicity of presentation, we set
$$\widehat{\mathfrak{E}}_k(t) =\int_{\mathcal{I}}  \bar\rho_0 |\pl_t^{k}w_t|^2 dx   +
  (1+t)^{-1}  \int_{\mathcal{I}}  {\bar\rho_0^\ga}   |\pl_t^k w_x|^2 dx.$$

Take $\pl_t^k$ onto \ef{eluerpert}, multiply the resulting equation by $  \pl_t^k w_t$,
and integrate the product  over  $\mathcal{I}$ to obtain
\be\label{516'-20}\begin{split}
\f{d}{dt} \overline{\mathfrak{E}}_k(t) +\int_{\mathcal{I}}  \bar\rho_0 |\pl_t^{k}w_t|^2 dx
  \le  C(\ga, M)(1+t)^{-2k-2} \lt( (\da+ \ea_0)   \mathcal{E}_{k}(t) \rt.\\
    \lt. + \lt(\da^{-1}+\ea_0\rt)  \lt(\widetilde{\mathcal{E}}_{0} + \sum_{1\le \iota\le k-1}\mathcal{E}_{\iota}\rt)(t) \rt)
      \end{split}\ee
for any positive number $\delta>0$,  where $\overline{\mathfrak{E}}_k$
satisfies  the following estimates:
 \begin{align*}
& \overline{\mathfrak{E}}_k \ge  C(\ga)^{-1} \widehat{\mathfrak{E}}_k -C(\ga, M) (1+t)^{-1-2k}\lt(\widetilde{\mathcal{E}}_{0} + \sum_{1\le \iota\le k-1}\mathcal{E}_{\iota} \rt),
    \\
& \overline{\mathfrak{E}}_k \le  C(\ga)\widehat{\mathfrak{E}}_k + C(\ga, M) (1+t)^{-1-2k}\lt(\widetilde{\mathcal{E}}_{0} + \sum_{1\le \iota\le k-1}\mathcal{E}_{\iota}\rt).
      \end{align*}
Take $\pl_t^k$ onto \ef{eluerpert}, multiply the resulting equation by $ \pl_t^k w$,
and integrate the product  over  $\mathcal{I}$ to obtain
\be\label{521-20}\begin{split}
&\f{d}{dt}  {\mathfrak{E}}_k(t)     + \widehat{\mathfrak{E}}_k (t)\le C(\ga, M)  (1+t)^{-1-2k} \lt(\widetilde{\mathcal{E}}_{0} + \sum_{1\le \iota\le k-1}\mathcal{E}_{\iota}\rt)(t),
      \end{split}\ee
where
\begin{align*}
&  \mathfrak{E}_k \ge C(\ga)^{-1} \lt( \int_{\mathcal{I}} \bar\rho_0 |\pl_t^k w|^2 dx    +\widehat{\mathfrak{E}}_k \rt)-C(\ga, M) (1+t)^{-1-2k} \lt(\widetilde{\mathcal{E}}_{0} + \sum_{1\le \iota\le k-1}\mathcal{E}_{\iota}\rt),
    \\
&  {\mathfrak{E}}_k
   \le C(\ga) \lt( \int_{\mathcal{I}} \bar\rho_0 |\pl_t^k w|^2 dx    +\widehat{\mathfrak{E}}_k \rt)-C(\ga, M) (1+t)^{-1-2k} \lt(\widetilde{\mathcal{E}}_{0} + \sum_{1\le \iota\le k-1}\mathcal{E}_{\iota}\rt).
\end{align*}
The derivation of estimates \ef{516'-20} and \ef{521-20} can be found in Lemma 3.6 of \cite{LZ}, which is based on \ef{decay-1d},  elliptic estimates \ef{20.el} and  $L^\infty$-estimates \ef{20lem31est}, so we omit details here. As in \cite{LZ}, we integrate the product of \ef{521-20} and $(1+t)^{m}$ over $[0,t]$ from $m=0$ to $m=2k$ step by step, and then integrate the product of \ef{516'-20} and $(1+t)^{1+2k}$ over $[0,t]$ to achieve
 \begin{align}
\mathcal{E}_k(t) +  \int_0^t  \mathcal{D}_{k}(s) ds
  \le C(\ga, M )\sum_{0\le \iota \le k}  \mathcal{E}_\iota (0).  \label{fal-1}
\end{align}
Based on \ef{12-4-j-20} and \ef{fal-1}, we integrate the product of \ef{521-20} and $\lt(\ln(1+t)\rt)^q(1+t)^{2k}$ over $[0,t]$ to give
 \begin{align}
&\lt(\ln(1+t)\rt)^q(1+t)^{2k} \int_{\mathcal{I}} \bar\rho_0 |\pl_t^k w|^2 dx
\notag\\
&+\int_0^t \lt(\ln(1+s)\rt)^q(1+s)^{2k-1} \int_{\mathcal{I}} \bar\rho_0^\ga |\pl_s^k w_x|^2 dx ds
  \le C(\ga, M ,q )  \mathcal{E} (0).  \label{fal-2}
\end{align}
We integrate the product of $\lt(\ln(1+t)\rt)^q(1+t)^{1+2k}$ and \ef{516'-20} over $[0,t]$, and use \ef{12-4-j-20}, \ef{fal-1} and \ef{fal-2} to obtain
\begin{align}
&\lt(\ln(1+t)\rt)^q(1+t)^{2k} \int_{\mathcal{I}} \lt( (1+t) \bar\rho_0 |\pl_t^k w_t|^2 + \bar\rho_0^\ga |\pl_t^k w_x|^2   \rt)dx \notag\\
&+\int_0^t \lt(\ln(1+s)\rt)^q(1+s)^{2k+1} \int_{\mathcal{I}} \bar\rho_0^\ga |\pl_s^k w_s|^2 dx ds
  \le C(\ga, M ,q )  \mathcal{E} (0).  \label{fal-3}
\end{align}
In fact, the derivations of \ef{fal-2} and \ef{fal-3} are similar to those of \ef{12.21} and \ef{date3.26}, respectively, in the three-dimensional spherically symmetric case, so we omit details here.

Using \ef{fal-1}-\ef{fal-3},  we prove that \ef{new12-4} holds for $j=k$. This finishes the proof of Lemma  \ref{lem-d3.8}.

\subsection{Proof of Theorem \ref{1d-mainthm1}}
Let $q>0$ be any fixed constant. It follows from \ef{wsv} and   $\bar\rho_0^{\ga-1}=A-Bx^2$   that
\begin{align*}
&\|w_x(\cdot, t) \|_{H^{1+(1+[\alpha]-\alpha)/{2}}(\mathcal{I})}^2\le C(\ga,M)\|w_x (\cdot, t)\|_{H^{\alpha+3+[\alpha],\ {3+[\alpha]}}(\mathcal{I})}^2
  \\
=& C(\ga,M) \sum_{0\le i \le 3+[\alpha]} \int_\mathcal{I} \bar\rho_0^{(\ga-1)(\alpha+3+[\alpha])} |\pl_x^i w_x|^2 dx
 \\
\le & C(\ga,M) \sum_{1\le i \le 4+[\alpha]} \int_\mathcal{I}  \bar\rho_0^{1+(i-1)(\ga-1)} |\pl_x^i w |^2 dx
 \\
\le & {C(\ga, M)} \sum_{1\le i\le 4+ [\alpha] } \mathcal{ E}_{0, i}(t),
\end{align*}
which, together with the embedding $H^{1}(\mathcal{I} )\hookrightarrow L^\infty(\mathcal{I} )$, and \ef{d3.12}, implies that
\begin{align}\label{Feb24-1}
\|w_x(\cdot, t) \|_{L^\infty(\mathcal{I})}^2\le C(\ga,M,2q) \lt( \ln(1+t)^{2q}+1 \rt)^{-1}\mathcal{E}(0)\notag\\
 \le C(\ga,M,q) \lt( \ln(1+t)^{q}+1 \rt)^{-2}\mathcal{E}(0), \ \  t\ge 0.
\end{align}
Similarly, we have for $t\ge 0$,
\begin{align}\label{fal-4}
 (1+t)^2  \|w_t(\cdot, t) \|_{L^\infty(\mathcal{I})}^2\le C(\ga,M,q) \lt( \ln(1+t)^{q}+1 \rt)^{-2}\mathcal{E}(0) .
\end{align}

Recall that  for $t\ge 0$,
\begin{subequations}\label{h-1d}
 \begin{align}
& 0\le h(t)\le C(\ga)  (1+t)^{-\frac{\ga}{\ga+1}}\ln(1+t),\label{h-1d-a}  \\
& |h_t|\le C(\ga) (1+t)^{-1-\frac{\ga}{\ga+1}}\ln(1+t). \label{h-1d-b}
 \end{align}\end{subequations}
 This is estimate (2.14) in \cite{LZ}, whose proof follows from the ODE analyses.
 Note that
\begin{align}
& \rho(\eta , t)-\bar\rho(\bar\eta , t)
 =\frac{\bar\rho_0 }{\eta_x }-\frac{\bar\rho_0 }{\bar\eta_x }
 = \frac{-\bar\rho_0 (  h+ w_x) }{(\bar\eta_x+h +w_x)\bar\eta_x }, \label{Feb28-4}\\
& u(\eta , t)-\bar u (\bar\eta , t)=w_t +xh_t,  \label{Feb28-6}
\end{align}
where \ef{Feb28-4} follows from   \ef{Feb25-1}, \ef{Feb25-3} and \ef{212}; and \ef{Feb28-6} from \ef{Feb26}, \ef{Feb28-5} and \ef{212}.
Then, \ef{1'-1d} follows from \ef{212},  \ef{h-1d-a} and \ef{Feb24-1}; and \ef{2'-1d} from \ef{fal-4} and \ef{h-1d-b}.

It follows from \ef{d3.10} and \ef{Feb24-1} that
\begin{align*}
\lt|w(x,t)-  \vartheta_0 \rt|
\le  C(\ga,M,q) \lt( \ln(1+t)^{q}+1 \rt)^{-1}\sqrt{\mathcal{E}(0)},
\end{align*}
which, together with  \ef{h-1d-a},  implies that  for $(x,t)\in \bar{\mathcal{I}}\times [0,\iy] $,
\begin{align}\label{Feb24-2}
 \lt|\eta(x,t)-\bar\eta(x,t) -\vartheta_0  \rt| \le C(\ga,M) (1+t)^{-\frac{\ga}{\ga+1}}\ln(1+t) \notag \\
  + C(\ga,M,q)\lt( \ln(1+t)^{q}+1 \rt)^{-1}\sqrt{\mathcal{E}(0)}.
\end{align}
\ef{Feb28-3} follows from \ef{vbs-1d},  \ef{212} and \ef{Feb24-2}.
Finally, \ef{4'-1d} follows from  $\tilde \eta=x \tilde \eta_x$,   $\ef{decay-1d-b}$,  \ef{20lem31est} and
$$\frac{d^k x_{\pm}(t)}{dt^k}=\pl_t^ k\tilde \eta \lt(\pm \sqrt{\frac{A_1}{B_1}}, t\rt)+ \pl_t^ k w \lt(\pm \sqrt{\frac{A_1}{B_1}}, t\rt), \  \ k=1, 2, 3.$$
This finishes the proof of Theorem \ref{1d-mainthm1}.

\section{The three-dimensional Spherically symmetric motions}\label{section3}
In this section, we consider the spherically symmetric free boundary problem in three dimensions, that is, we seek solutions with symmetry to  problem \ef{m2.1} of the form:
$$\Omega(t) =B_{R(t)}(  0 ), \ \  \rho(x, t)=\rho(r, t), \ \ u(x, t)= (x/r) u(r, t),  $$
where  $r=|x|$ and $B_{R(t)}(  0 )$ is the ball centered at the origin with the radius $R(t)$.
Then  problem \ef{m2.1} reduces to
\begin{subequations}\label{equation1}
\begin{align}
&    (r^2\rho)_t+ (r^2\rho u)_r=0    & {\rm in } & \  \  \lt(0, \  R(t)\rt), \label{equation1a}\\
&\rho( u_t +u  u_r)+ p_r=- \rho u  & {\rm in } & \  \  \lt(0, \  R(t)\rt),  \label{equation1b}\\
&\rho > 0  & {\rm in} &   \  \  \lt[0, \  R(t)\rt),\label{equation1c}\\
& \rho\lt(R(t),t\rt) =0,      \ \ u(0,t)=0,     \label{equation1d}\\
&\dot R(t)=u(R(t), t) \ \ {\rm with} \ \ R(0)=R_0,   \label{equation1e}\\
& (\rho,u)(r,t=0)=\lt(\rho_0, u_0\rt)(r) & {\rm  on } & \  \  \lt(0, \  R_0\rt), \label{equation1f}
\end{align}
\end{subequations}
so that $R(t)$ is the  radius of the domain occupied by the gas at time $t$, and $r=R(t)$  represents  the vacuum boundary which issues from $r=R_0$ and moves with the fluid velocity.
In the setting, the initial domain is taken to be the  ball  $B_{R_0}(0)$ and the condition \ef{initial density} satisfied by the initial density reduces to
\begin{subequations}\label{March2}
\begin{align}
& \rho_0(r)>0 \ \ {\rm for} \ \  0\le r<R_0 , \ \ \rho_0(R_0)=0, \ \ 4\pi \int_0^{R_0} r^2 \rho_0(r)dr =M, \\
&  -\iy<   \lt(\rho_0^{\ga-1}\rt)_r <0 \  \ {\rm at} \  \ r=R_0.
 \end{align}
\end{subequations}

As in \cite{HZ},
we transform problem \ef{equation1} into Lagrangian variables to fix the boundary, and make  the  initial ball  of the Barenblatt solution, $\bar\Omega(0)=B_{\sqrt{A_3/B_3}}(0)$,
as the reference domain. For simplicity of presentation,  set
$$ \mathcal{I}   =\lt(0, \ \sqrt{A_3 /B_3} \rt). $$
For $r\in \mathcal{I}$, we define the Lagrangian variable  $\eta(r, t)$ by
\be\label{haz2.2}
  \eta_t(r, t)= u(\eta(r, t), t) \ \ {\rm for} \  \ t>0, \  \ {\rm and} \ \  \eta(r, 0)=\eta_0(r) .
\ee
Then, it follows from \ef{equation1a} that $\pl_t(\rho(\eta(r, t), t) \eta^2(r,t) \eta_r(r,t))=0$, which implies that
$$\rho(\eta(r, t), t) \eta^2(r,t) \eta_r(r,t)=\rho_0(\eta_0(r)) \eta_0^2(r) \eta_{0r}(r).$$
 If we choose $\eta_0$ such that $\rho_0(\eta_0(r)) \eta_0^2(r) \eta_{0r}(r)=r^2 \bar\rho(r,0)$, where $\bar\rho(r,0)$ is the initial density of the Barenblatt solution, then we have for $r\in \mathcal{I}$,
\be\label{newlagrangiandensity}
\rho(\eta(r, t), t)
 =\frac{r^2 \bar\rho_0(r)}{\eta^2(r,t) \eta_r(r,t)}, \ \ {\rm where} \ \
 \bar\rho_0(r) = \bar\rho(r,0)=(A_3-B_3r^2)^{\frac{1}{\ga-1}}.\ee
Such an $\eta_0$ exists, for instance, we may define
$\eta_0: \mathcal{I} \to \lt(0, R_0\rt)  $ by
\bee
\int_0^{\eta_0(r)} r^2 \rho_0(r) dr = \int_0^{r} r^2 \bar\rho_0(r)dr \ \  {\rm for} \ \
r\in  \mathcal{I}.
\eee
Due to \ef{Feb25-1}, \ef{Feb25-4} and \ef{March2}, $\eta_0$ is a well-defined diffeomorphism.
So, the vacuum free boundary problem \ef{equation1} is reduced to the following initial-boundary value problem with the fixed boundary:
\begin{subequations}\label{419} \begin{align}
& \bar\rho_0\eta_{tt} + \bar\rho_0 \eta_t  + \lt( \frac{\eta}{r}\rt)^2  \left[    \lt(\frac{r^2}{\eta^2}\frac{\bar\rho_0}{ \eta_r}\rt)^\ga    \right]_r =0   \ \  &{\rm in}& \  \    \mathcal{I}\times (0, \iy), \label{419a}\\
&  \eta(0, t)=0,    &  {\rm on}& \  \ (0,\iy),\\
& (\eta, \eta_t) = \lt(\eta_0, u_0(\eta_0)\rt)  &  {\rm on}& \  \    \mathcal{I}\times (0, \iy).
\end{align}
\end{subequations}
In the setting, the radius of the moving vacuum  boundary for problem \ef{equation1} is given by
\be\label{vbs}  R(t)= \eta\lt(\sqrt{A_3/B_3}, \ t\rt) \ \  {\rm for} \ \  t\ge 0.\ee

\subsection{Main results}
Define the Lagrangian variable $\bar\eta(r, t)$ for the Barenblatt flow in $\bar {\mathcal{I}}$ by
\be\label{haz2.7}
 \pl_t \bar \eta(r, t)= \bar u(\bar \eta (r, t), t)=\f{  \bar\eta(r,t)}{(3\ga-1)(1+ t)}  \ \ {\rm for} \  \ t>0 \  \ {\rm and} \ \  \bar \eta(r, 0)=r,
\ee
then
\be\label{bareta}
\bar \eta(r,t)= r \lt(1 +  t \rt)^{\frac{1}{3\ga-1}}   \ \  {\rm for} \ \  (r,t)\in \bar{\mathcal{I}}\times [0,\iy).
\ee
Since $\bar\eta$ does not solve \ef{419a} exactly,
a correction $h(t)$ was introduced in \cite{HZ} which is a solution of the following initial value problem of ordinary differential equations:
\begin{align*}
& h_{tt} + h_t - (3\ga-1)^{-1}  (\bar \eta_r+h)^{2-3\ga}   + \bar\eta_{rtt}  +\bar\eta_{rt}    =0, \ \ t >0,  \\
& h(t=0)=h_t(t=0)=0.
\end{align*}
The new ansatz is then given by
$
\tilde \eta(r, t)=\bar \eta (r,t)  + r h(t)
$,
so that
 \be\label{equeta} \begin{split}
&  \bar\rho_0  \tilde\eta_{tt}   +  \bar\rho_0  \tilde\eta_t + \lt(\frac{\tilde\eta}{r}\rt)^2 \left[    \lt(\frac{r^2}{\tilde\eta^2}\frac{\bar\rho_0}{ \tilde\eta_r}\rt)^\ga    \right]_r =0    \ \  {\rm in} \  \  {\mathcal{I}}\times (0,\iy).
\end{split}
\ee
It was proved in \cite{HZ} that $\tilde\eta$ behaves similarly to $\bar\eta$, see \ef{decay} and \ef{decayforh} for details. We write equation \ef{419a} in a perturbation form around the Barenblatt solution. Let
$\zeta(r,t)= r^{-1}\eta(r,t)- r^{-1} {\tilde\eta(r,t)} $,
then it follows from \ef{419a} and \ef{equeta} that
\be\label{pertb}\begin{split}
   &r \bar\rho_0 \zeta_{tt}   +
   r  \bar\rho_0  \zeta_t
   +   \lt(\tilde\eta_r+ \zeta\rt)^2\left[ \bar\rho_0^\ga     \lt(\tilde\eta_r+\zeta\rt)^{-2\ga}\lt(\tilde\eta_r+\zeta +r\zeta_r\rt)^{-\ga}   \right]_r -    \tilde\eta_r^{2-3\ga}  \lt(   {\bar\rho_0}  ^\ga    \right)_r
   =0 .
\end{split}\ee
We denote
$ \alpha = (\ga-1)^{-1}$, and set that for
$j=0,\cdots, 4 +[\alpha]$ and  $i=0,\cdots, 4 +[\alpha]-j$,
\bee\label{}\begin{split}
&\mathcal{ E}_{j}(t)  = (1+ t)^{2j} \int_\mathcal{I}   \lt[r^4 \bar\rho_0  \lt(\pl_t^j \zeta\rt)^2 +r^2 \bar\rho_0^\ga  \lt|\pl_t^j \lt( \zeta, r\zeta_r \rt)  \rt|^2  + (1+ t)  r^4 \bar\rho_0    \lt(\pl_t^{j} \zeta_t\rt)^2 \rt]   dr,\\
&\mathcal{ E}_{j, i}(t) =  (1+  t)^{2j}  \int_\mathcal{I} \lt[r^2 \bar\rho_0^{1+(i-1)(\ga-1)}  \lt(\pl_t^j \pl_r^i \zeta\rt)^2 + r^4 \bar\rho_0^{\ga+i(\ga-1)}  \lt(\pl_t^j \pl_r^{i} \zeta_r\rt)^2 \rt] dr  .
\end{split}\eee
 The higher-order norm is defined by
\be\label{d3.29-2}
  \mathcal{E}(t) =  \sum_{0\le j\le 4+[\alpha]} \lt(\mathcal{ E}_{j}(t) + \sum_{1\le i \le 4+[\al]-j} \mathcal{ E}_{j, i}(t)  \rt).
  \ee

 The main results in this section are stated as follows.
{\begin{thm}\label{thm1} There exists a  constant $\epsilon_0>0$ such that if
$\mathcal{E}(0)\le \epsilon_0^2$,
then   problem \eqref{419}  admits a global unique smooth solution  in $\mathcal{I}\times[0, \iy)$ satisfying that for any fixed constant $q>0$,
\be\label{12.16}
\lt(\lt(\ln(1+t)\rt)^q +1\rt)\mathcal{E}(t)\le C(\ga, M, q)\mathcal{E}(0), \ \  \ \  t\ge 0.
\ee
\end{thm}

With Theorem \ref{thm1}, we have the following theorem for solutions to the original vacuum  free boundary problem \ef{equation1} concerning  the time-asymptotics of the vacuum boundary and better convergence rates of the density and velocity  than those in \cite{HZ}.

\begin{thm}\label{thm2} There exists a  constant $\epsilon_0>0$ such that if
$\mathcal{E}(0)\le \epsilon_0^2$,
then  problem
\ef{equation1}  admits a global unique smooth solution $\lt(\rho, u,  R(t)\rt)$ for $t\in[0,\iy)$. Let $q>0$ be any fixed constant, and set
$$\mathscr{H}(t)=\lt(\lt(\ln(1+t)\rt)^q +1\rt)^{-1}\sqrt{\mathcal{E}(0)}+(1+t)^{-\frac{3\ga-2}{3\ga-1}}\ln(1+t), $$
then it holds that for all $r \in \mathcal{I}$ and $t\ge 0 $,
\begin{subequations}\label{date3.28}\begin{align}
 &\lt|\rho\lt(\eta(r, t),t\rt)-\bar\rho\lt(\bar\eta(r, t), t\rt)\rt|
 \le  C(\ga, M,q)\lt(A_3-B_3 r ^2\rt)^{\frac{1}{\ga-1}}(1+t)^{-\frac{4}{3\ga-1}} \mathscr{H}(t),
 \label{thm2est1}\\
 & \lt|u\lt(\eta(r, t),t\rt)-\bar u\lt(\bar\eta(r, t), t\rt)\rt| \le  C(\ga, M,q) r (1+t)^{-1} \mathscr{H}(t);
\label{thm2est2}
\end{align}\end{subequations}
 and for all $t\ge 0$,
\begin{subequations}\begin{align}
&\lt| R(t) - \bar R (t) \rt|\le C(\ga, M,q) \mathscr{H}(t) ,\label{thm2est3}\\
& \lt|\frac{d^k R(t)}{dt^k}\rt|\le C(\ga, M)(1+t)^{\frac{1}{3\ga-1}-k} , \ \   k=1, 2 ,3 ,
\label{thm2est4}
\end{align}\end{subequations}
where $\bar R(t)=\sqrt{A_3/B_3}(1+t)^{1/(3\ga-1)}$ is the radius of the vacuum boundary of the Barenblatt solution.
\end{thm}

\begin{rmk} In \cite{HZ}, it was proved that the estimates in \ef{date3.28} and \ef{thm2est3} hold with
$$\mathscr{H}(t)=\sqrt{\mathcal{E}(0)}+(1+t)^{-\frac{3\ga-2}{3\ga-1}}\ln(1+t). $$
So,  the estimates in  \ef{date3.28} improve the convergence rates obtained in \cite{HZ}; and estimate \ef{thm2est3} gives the convergence of $R(t)$  to  $\bar R(t)$, in addition to the same expanding rate of $R(t)$ and $\bar R(t)$ obtained in \cite{HZ}.
\end{rmk}

\subsection{Proof of Theorem \ref{thm1}}
It follows from Theorem 2.1 in \cite{HZ} that there exists a small constant $\epsilon_0 >0$ such that if $\mathcal{E}(0)\le \epsilon_0^2$, then problem \eqref{419}  admits a global unique smooth solution $\zeta(r,t)$ in $\mathcal{I}\times[0, \iy)$ satisfying
\be\label{10.30}
\mathcal{E}(t)\le C(\ga, M)\mathcal{E}(0)   \  \ {\rm for} \ \    t\ge 0.
\ee
So, it suffices to prove  the following estimates:
\begin{lem}\label{lem1218} Let $q$ be any fixed positive integer, then we have for all $t\ge 0$,
  \begin{align}
\lt(\lt( \ln(1+t) \rt)^q +1 \rt)\mathcal{E}_j(t) +  \int_0^t \lt( \lt( \ln(1+s) \rt)^q +1 \rt) \mathcal{D}_{j}(s) ds  \notag \\
  \le C(\ga, M,q)  \sum_{0\le \iota \le j} \mathcal{E}_\iota (0), \ \  \ \  j=0 ,1 ,2, \cdots,  4 +[\alpha],  \label{12-4}
\end{align}
 where
\begin{align*}
\mathcal{D}_{j}(t)   = (1+ t)^{2j-1} \int_\mathcal{I}    r^2 \bar\rho_0^\ga  \lt|\pl_t^j \lt( \zeta, r\zeta_r \rt)  \rt|^2  dr + (1+ t)^{2j+1} \int_\mathcal{I} r^4 \bar\rho_0    \lt(\pl_t^{j} \zeta_t\rt)^2    dr.
\end{align*}
\end{lem}
Once this lemma is proved, estimate \ef{12.16} will be an immediate conclusion of the fact that $ ( \ln(1+t)  )^q \le C(q)  ( ( \ln(1+t)  )^{[q]+1} +1 )$ for any positive $q\notin \mathbb{Z}$ and the following  elliptic estimates:
\begin{align}\label{prop1est}
   \mathcal{E}_{j, i}(t)\le C(\ga, M)  \sum_{0\le \iota \le i+j} \mathcal{E}_\iota(t) \ \ {\rm when} \ \  j\ge 0, \ \ i\ge 1, \ \  i+j\le 4 +[\alpha].
\end{align}
This is Proposition 3.1 in \cite{HZ}, whose proof is based on the following $L^\iy$-estimates:
\begin{subequations}\label{lem34est}\begin{align}
&S(t)\le C(\ga,M)\mathcal{E}(t),\label{lem34est-a}\\
&S(t)\le C(\ga, M)\mathcal{E}(0)\le C(\ga, M) \epsilon_0^2, \label{lem34est-b}
\end{align}\end{subequations}
where
\begin{align*}
S(t)=  &  \sum_{j=0,1,2} (1+  t)^{2j} \lt \|\pl_t^j \zeta (\cdot,t) \rt \|_{L^\iy}^2   +  \sum_{j=0,1} (1+  t)^{2j} \lt \|\pl_t^j \pl_r \zeta  (\cdot,t) \rt\|_{L^\iy}^2
    \notag\\
    & +
  \sum_{
  i+j\le 2+[\al],\   2i+j \ge 3} (1+  t)^{2j}      \lt \|  \bar\rho_0^{\f{(\ga-1)(2i+j-3)}{2}}\pl_t^j \pl_r^i \zeta(\cdot,t) \rt\|_{L^\iy}^2 \notag\\
    & +
  \sum_{
  i+j=3+[\al]} (1+  t)^{2j}   \lt \| r \bar\rho_0^{\f{(\ga-1)(2i+j-3)}{2}}\pl_t^j \pl_r^i \zeta  (\cdot,t) \rt\|_{L^\iy}^2
   \notag \\
 & + \sum_{
  i+j=4+[\al]} (1+  t)^{2j}   \lt \| r^2 \bar\rho_0^{\f{(\ga-1)(2i+j-3)}{2}}\pl_t^j \pl_r^i \zeta (\cdot,t) \rt\|_{L^\iy}^2 .
\end{align*}
Indeed, \ef{lem34est-a} follows from Lemma 3.7 in \cite{HZ}, and \ef{lem34est-b}  from \ef{lem34est-a} and \ef{10.30}.

{\em Proof of Lemma \ref{lem1218}}. The proof consists of two steps. In {\em Step 1}, we prove \ef{12-4} for $j=0$, the zeroth order estimate. Based on this, we prove \ef{12-4} for $j\ge 1$, the higher order estimates in {\em Step 2},
since the time decay rate of the zeroth order norm $\mathcal{E}_0(t)$ determines those of the higher order norms $\mathcal{E}_j(t)$ ($j\ge 1$).

We first recall that  for all $t\ge 0$,
\begin{subequations}\label{decay}\begin{align}
&\lt(1 +   t \rt)^{\frac{1}{3\ga-1}} \le \tilde \eta_{r}(t) \le C(\ga) \lt(1 +   t \rt)^{\frac{1}{3\ga-1}},  \ \   \ \   \tilde\eta_{rt} \ge 0 , \label{decaya} \\
&\lt|\f{d^k\tilde \eta_{r}(t)}{dt^k}\rt| \le C(\ga, k)\lt(1 +   t \rt)^{\frac{1}{3\ga-1}-k},   \ \ k=1, 2, 3, \cdots.  \label{decayb}
 \end{align}\end{subequations}
Indeed, this is estimate (2.13) in \cite{HZ}, whose proof  follows from the ODE analyses.
For simplicity of presentation, we set
$$\widehat{\mathfrak{E}}_j(t)=(1+ t)^{-1} \int_\mathcal{I}    r^2 \bar\rho_0^\ga  \lt|\pl_t^j \lt( \zeta, r\zeta_r \rt)  \rt|^2  dr +  \int_\mathcal{I} r^4 \bar\rho_0    \lt(\pl_t^{j} \zeta_t\rt)^2    dr,$$
then
\begin{subequations}\label{12.19}\begin{align}
& \mathcal{E}_j(t)=(1+t)^{2j+1}\widehat{\mathfrak{E}}_j(t)+ (1+t)^{2j}\int_\mathcal{I} r^4 \bar\rho_0    \lt(\pl_t^{j} \zeta \rt)^2    dr,
\label{12.19a}\\
& \mathcal{D}_j(t)  \ge (1+t)^{2j}\widehat{\mathfrak{E}}_j(t) .\label{12.19b}
\end{align}\end{subequations}

{\em Step 1}. In this step, we prove \ef{12-4} for $j=0$, that is,
  \begin{align}
\lt(\lt( \ln(1+t) \rt)^q +1 \rt)\mathcal{E}_0 (t) +  \int_0^t \lt( \lt( \ln(1+s) \rt)^q +1 \rt) \mathcal{D}_{0}(s) ds    \le  C(\ga, M, q)  \mathcal{E}_0 (0) . \label{12-4-1}
\end{align}

Multiply \ef{pertb} by $r^3 \zeta_t$ and integrate the resulting equation over $\mathcal{I}$ to get
\begin{align}\label{4.2}
 \f{d}{dt}\int_\mathcal{I} \lt( \f{1}{2}r^4\bar\rho_0 \zeta_{t}^2 + r^2 \bar\rho_0^\ga L_1 \rt) dr   +
  \int_\mathcal{I} r^4 \bar\rho_0  \zeta_t^2 dr =-
  \tilde\eta_{rt}\int_\mathcal{I} r^2\bar\rho_0^\ga    {F}   dr
   ,
\end{align}
where
\begin{align*}
  L_1= & ({ \ga-1})^{-1} \lt(\lt(\tilde\eta_r+\zeta\rt)^{2-2\ga }\lt(\tilde\eta_r+\zeta +r\zeta_r\rt)^{1-\ga }  -\tilde\eta_r^{3-3\ga } \rt. \\
& \lt.  + (\ga-1) \tilde\eta_r^{2-3\ga} \lt( 3\zeta + r\zeta_r\rt) \rt)
,\\
 {F} = & 2\lt(\tilde\eta_r+\zeta\rt)^{1-2\ga }\lt(\tilde\eta_r+\zeta +r\zeta_r\rt)^{1-\ga }  +  \lt(\tilde\eta_r+\zeta\rt)^{2-2\ga }\lt(\tilde\eta_r+\zeta +r\zeta_r\rt)^{-\ga} \\
  & -3 \tilde\eta_r^{2-3\ga}
   -(2-3\ga) \tilde\eta_r^{1-3\ga}  \lt( 3\zeta + r\zeta_r\rt).
\end{align*}
It follows from the Taylor expansion, \ef{lem34est-b} and \ef{decaya}   that
\begin{align*}
& \lt| L_1 - 2^{-1} \tilde\eta_r^{1-3\ga}
 \mathscr{Q} \rt| \le C(\ga)
     \tilde\eta_r^{ -3\ga} \lt(|\zeta|^3+|r\zeta_{r}|^3\rt)  ,\\
&  {F}  \ge  2^{-1}(3\ga-1) \tilde\eta_r^{ -3\ga} \mathscr{Q}  - C(\ga)  \tilde\eta_r^{ -3\ga-1} \lt(|\zeta|^3+|r\zeta_{r}|^3\rt)  ,
 \end{align*}
where
\begin{align*}
\mathscr{Q}= 3\lt( 3 \ga-2\rt) \zeta^2 + 2(3\ga-2) \zeta r\zeta_r + \ga (r\zeta_r)^2 .
 \end{align*}
Clearly, $\mathscr{Q}$ is equivalent to $\zeta^2+ \lt|r\zeta_r\rt|^2$, that is,
\begin{align}\label{10.25-6}
(C(\ga))^{-1} \lt(  \zeta^2+ \lt|r\zeta_r\rt|^2 \rt)  \le \mathscr{Q} \le C(\ga)  \lt(  \zeta^2+ \lt|r\zeta_r\rt|^2 \rt),
 \end{align}
 which implies that
\begin{subequations}\label{10.25-1}\begin{align}
&(C(\ga))^{-1}    \lt(  \zeta^2+ \lt|r\zeta_r\rt|^2 \rt)  \le   (1+t) L_1  \le  C(\ga)   \lt(  \zeta^2+ \lt|r\zeta_r\rt|^2 \rt) ,\\
& {F}\ge C(\ga)^{-1} (1+t)^{-1-\frac{1}{3\ga-1}} \lt(  \zeta^2+ \lt|r\zeta_r\rt|^2 \rt)  \ge 0 .
 \end{align}\end{subequations}
Integrate \ef{4.2} over $[0,t]$, and use \ef{decaya} and \ef{10.25-1} to  arrive at
\begin{align}
   \widehat{\mathfrak{E}}_0(t)   +
  \int_0^t \int_\mathcal{I} r^4 \bar\rho_0  \zeta_s^2 drds
  \le C(\ga) \widehat{\mathfrak{E}}_0(0) . \label{4.4}
\end{align}

Multiply \ef{pertb} by $r^3 \zeta$ and integrate the resulting equation over $\mathcal{I}$ to obtain
\be\label{4.5}\begin{split}
  \frac{1}{2}\f{d}{dt}\int_\mathcal{I}  r^4\bar\rho_0 \lt( \zeta^2 + 2\zeta \zeta_t \rt) dr
   + \int_\mathcal{I} r^2 \bar\rho_0^\ga  L_2 dr
   = \int_\mathcal{I} r^4 \bar\rho_0  \zeta_t^2 dr,
\end{split}\ee
where
\begin{align*}
 L_2
    =& \lt[3\tilde\eta_r^{2-3\ga}-2 \lt(\tilde\eta_r+\zeta\rt)^{1-2\ga }\lt(\tilde\eta_r+\zeta +r\zeta_r\rt)^{1-\ga }
 \rt.\\
  &\lt. -\lt(\tilde\eta_r+\zeta\rt)^{2-2\ga }\lt(\tilde\eta_r+\zeta +r\zeta_r\rt)^{-\ga}\rt]\zeta  \\
  &+ \lt[\tilde\eta_r^{2-3\ga}-\lt(\tilde\eta_r+\zeta\rt)^{2-2\ga }\lt(\tilde\eta_r+\zeta +r\zeta_r\rt)^{-\ga}\rt]r\zeta_r   .
\end{align*}
Again, we use the Taylor expansion, \ef{lem34est-b}, \ef{decaya} and \ef{10.25-6}   to  give
\begin{align}
 L_2
\ge \tilde\eta_r^{1-3\ga} \mathscr{Q} -C(\ga)  \tilde\eta_r^{ -3\ga} \lt(|\zeta|^3+|r\zeta_{r}|^3\rt)
    \ge C(\ga)^{-1} (1+t)^{-1} \lt( \zeta^2+ \lt(r\zeta_r\rt)^2 \rt) . \label{12.10-2}
\end{align}
Then, we integrate \ef{4.5} over $[0,t]$,  and use the Cauchy inequality and \ef{4.4} to achieve
\be\label{4.6}\begin{split}
\int_\mathcal{I}   \lt( r^4\bar\rho_0 \zeta ^2\rt)  (r,t) dr   +
  \int_0^t  (1+s)^{-1} \int_\mathcal{I}  r^2 \bar\rho_0^\ga \lt(\zeta^2+ \lt(r\zeta_r\rt)^2\rt)  drds
  \le C(\ga) \mathcal{E}_0(0).
\end{split}\ee

Integrate the product of \ef{4.2} and (1+t) over $[0,t]$, and use \ef{decaya}, \ef{10.25-1},  \ef{4.4} and \ef{4.6} to get
\begin{align*}
(1+t) \widehat{\mathfrak{E}}_0(t)
  +
 \int_0^t (1+s) \int_\mathcal{I} r^4 \bar\rho_0  \zeta_s^2 drds \notag \\
 \le C(\ga) \widehat{\mathfrak{E}}_0(0) + \int_0^t \widehat{\mathfrak{E}}_0(s) ds
  \le C(\ga) \mathcal{E}_0(0),
\end{align*}
which, together with \ef{12.19a} and \ef{4.6}, implies
  \begin{align}
 \mathcal{E}_0(t) +  \int_0^t   \mathcal{D}_{0}(s) ds   \le  C( \ga)  \mathcal{E}_0 (0) . \label{12.4-1}
\end{align}

Let $m$ be any given positive integer, then we integrate the product of  $\lt(\ln(1+t)\rt)^m $ and \ef{4.5} over $[0,t]$  to obtain
\begin{align}\label{12.10-1}
&  \lt(\ln(1+t)\rt)^m \int_\mathcal{I}   \lt( r^4\bar\rho_0 \zeta ^2\rt)  (r,t) dr  \notag\\
 &  +
  \int_0^t \lt(\ln(1+s)\rt)^m  (1+s)^{-1}  \int_\mathcal{I}  r^2 \bar\rho_0^\ga \lt(\zeta^2+ \lt(r\zeta_r\rt)^2\rt)  drds \notag\\
&  \le C(\ga) \lt\{    \mathcal{E}_0(0) + \lt(\ln(1+t)\rt)^m \int_\mathcal{I}   \lt( r^4\bar\rho_0 \zeta_t^2\rt)  (r,t) dr \rt. \notag\\
 & + m \int_0^t \lt(\ln(1+s)\rt)^{m-1}(1+s)^{-1}  \int_\mathcal{I} r^4\bar\rho_0 ( \zeta ^2 + \zeta_s^2) dr ds \notag\\
 &\lt. +
  \int_0^t  \lt(\ln(1+s)\rt)^m  \int_\mathcal{I} r^4 \bar\rho_0  \zeta_s^2 drds  \rt\} \notag\\
 &  \le   C( \ga, m)  \mathcal{E}_0(0) +  C(\ga, M)  m \notag \\
  & \times \int_0^t \lt(\ln(1+s)\rt)^{m-1}(1+s)^{-1} \int_\mathcal{I}  r^2 \bar\rho_0^\ga \lt(\zeta^2+ \lt(r\zeta_r\rt)^2\rt) dr ds  ,
\end{align}
where the first inequality follows from the Cauchy inequality and \ef{12.10-2}, and  the second  from \ef{12.4-1} and the following estimate:
\begin{align}
\int_\mathcal{I}  r^4\bar\rho_0 \zeta ^2 dr \le C(\ga, M) \int_\mathcal{I}  r^2 \bar\rho_0^\ga \lt(\zeta^2+ \lt(r\zeta_r\rt)^2\rt)  dr. \label{12-19}
\end{align}
Clearly, it follows from  \ef{12.10-1} and \ef{4.6} inductively that
\begin{align}\label{12.9}
& \int_0^t  \lt(\ln(1+s) \rt)^q    (1+s)^{-1} \int_\mathcal{I}  r^2 \bar\rho_0^\ga \lt(\zeta^2+ \lt(r\zeta_r\rt)^2\rt)  drds \notag\\
  & +\lt(\ln(1+t) \rt)^q \int_\mathcal{I}   \lt( r^4\bar\rho_0 \zeta ^2\rt)  (r,t) dr \le C( \ga,M, q)  \mathcal{E}_0(0).
\end{align}
It remains to prove \ef{12-19}.  Since $\bar\rho_0$ has the positive minimum and maximum on  $[0, \bar a/2]$ with $\bar a = \sqrt{A_3/B_3}$,  then
 \begin{align*}
  \int_0^{\bar a/2} r^4\bar\rho_0 \zeta ^2 dr  \le C(\ga, M) \int_0^{\bar a/2} r^2 \zeta ^2  dr
  \le C(\ga, M)  \int_0^{\bar a/2} r^2 \bar\rho_0^\ga \zeta ^2  dr .
\end{align*}
Near the boundary, it follows from the Hardy inequality \ef{hardy2019}  and the equivalence of $\bar\rho_0^{\ga-1}$ and $\bar a-r$ that
\begin{align*}
  \int_{\bar a/2}^{\bar a} r^4\bar\rho_0 \zeta ^2 dr\le C(\ga, M) \int_{\bar a/2}^{\bar a} \bar\rho_0 \zeta^2 dr \le C(\ga, M)   \int_{\bar a/2}^{\bar a} \bar\rho_0^{1+2(\ga-1)}
(\zeta ^2 +\zeta_r^2 ) dr \\
 \le  C(\ga, M)  \int_{\bar a/2}^{\bar a} \bar\rho_0^{\ga}
(\zeta ^2 +\zeta_r^2 ) dr  \le  C(\ga, M)  \int_{\bar a/2}^{\bar a}  r^2 \bar\rho_0^\ga \lt(\zeta^2+ \lt(r\zeta_r\rt)^2\rt)  dr.
\end{align*}
This finishes the proof of \ef{12-19}.
Finally, we integrate the product of $\lt(\ln(1+t)\rt)^q (1+t)$  and  \ef{4.2} over $[0,t]$ to achieve
\begin{align*}
&\lt(\ln(1+t)\rt)^q  (1+t)   \widehat{\mathfrak{E}}_0(t)
    +
 \int_0^t  \lt(\ln(1+s)\rt)^q (1+s) \int r^4 \bar\rho_0  \zeta_s^2 drds \notag\\
&  \le  C(\ga) \widehat{\mathfrak{E}}_0(0) +
 C(\ga) \int_0^t \lt(q\lt(\ln(1+s) \rt)^{q-1} +\lt(\ln(1+s) \rt)^{q}\rt) \widehat{\mathfrak{E}}_0(s) ds \notag\\
&  \le   C( \ga,M, q)  \mathcal{E}_0(0) ,
\end{align*}
where the first inequality follows from \ef{10.25-1}, and the second from \ef{12.4-1} and \ef{12.9}.
This, together with \ef{12.9} and \ef{12.4-1}, proves \ef{12-4-1}.

{\em Step 2}.  In this step, we use the mathematical induction to prove \ef{12-4} for $j \ge 1$. It follows from \ef{12-4-1} that \ef{12-4} holds for $j=0$. For $1\le k\le 4+ [\alpha]$, we make the hypothesis that \ef{12-4} holds for $j=0, 1, \cdots, k-1$, that is,
 \begin{align}
\lt(\lt( \ln(1+t) \rt)^q +1 \rt)\mathcal{E}_j(t) +  \int_0^t \lt( \lt( \ln(1+s) \rt)^q +1 \rt) \mathcal{D}_{j}(s) ds  \notag \\
  \le C(\ga, M,q )  \sum_{0\le \iota \le j} \mathcal{E}_\iota (0), \ \  \ \  j=0 ,1 ,2, \cdots,  k-1.  \label{12-4-j}
\end{align}
It is enough to prove \ef{12-4} holds for $j=k$ under the hypothesis \ef{12-4-j}.

Take $\pl_t^k$ onto \ef{pertb}, multiply the resulting equation by $r^3 \pl_t^k \zeta_t$,
and integrate the product  over  $\mathcal{I}$ to obtain
\be\label{3-62}\begin{split}
   & \frac{d}{dt}\overline{\mathfrak{E}}_k (t) + \int_\mathcal{I}  r^4  \bar\rho_0  \lt(\pl_t^k \zeta_t\rt)^2 dr
  \\     \le &   C(\ga, M)(1+t)^{  -2-2k} \lt(   (\ea_0+\da) \mathcal{E}_k(t) +   \lt(\ea_0+\da^{-1}\rt)   \sum_{0\le \iota \le k-1}  \mathcal{E}_\iota(t) \rt)
\end{split}\ee
for any positive number $\da>0$, where $\overline{\mathfrak{E}}_k$ satisfies the following estimates:
\begin{align*}
&    \overline{\mathfrak{E}}_k \ge C(\ga)^{-1} \widehat{\mathfrak{E}}_k - C(\ga, M) (1+t)^{-1-2k} \sum_{0\le \iota\le k-1}  \mathcal{E}_\iota(t) , \\
&   \overline{\mathfrak{E}}_k \le C(\ga) \widehat{\mathfrak{E}}_k  + C(\ga, M) (1+t)^{-1-2k} \sum_{0\le \iota\le k-1}  \mathcal{E}_\iota(t)  .
\end{align*}
Take $\pl_t^k$ onto \ef{pertb}, multiply the resulting equation by $r^3 \pl_t^k \zeta$,
 integrate the product  over  $\mathcal{I}$, and use \ef{3-62} to get
\be\label{3-73}\begin{split}
      \frac{d}{dt} \mathfrak{E}_k(t) +   \widehat{\mathfrak{E}}_k(t)
  \le C(\ga, M)    \sum_{0\le \iota\le k-1}  (1+ t)^{ -2k+2\iota} \widehat{\mathfrak{E}}_\iota(t) ,
\end{split}\ee
where
\begin{align*}
  \mathfrak{E}_k  \ge & C(\ga)^{-1}\lt( \int_\mathcal{I} r^4  \bar\rho_0 \lt| \pl_t^k  \zeta    \rt|^2 dr +  \widehat{\mathfrak{E}}_k(t) \rt)  - C(\ga, M) (1+t)^{-1-2k} \sum_{0\le \iota\le k-1}  \mathcal{E}_\iota(t), \\
 \mathfrak{E}_k  \le & C(\ga) \lt( \int_\mathcal{I} r^4  \bar\rho_0 \lt| \pl_t^k  \zeta    \rt|^2 dr + \widehat{\mathfrak{E}}_k(t) \rt) + C(\ga, M)  (1+t)^{-1-2k} \sum_{0\le \iota\le k-1}  \mathcal{E}_\iota(t)
.
\end{align*}
Indeed, the derivation of estimates \ef{3-62} and \ef{3-73} can be found in Lemma 3.6 of \cite{HZ}, which is based on \ef{decay},  elliptic estimates \ef{prop1est} and  $L^\infty$-estimates \ef{lem34est}, so we omit details here. As in \cite{HZ}, we integrate the product of \ef{3-73} and $(1+t)^{m}$ over $[0,t]$ from $m=0$ to $m=2k$ step by step, and then integrate the product of \ef{3-62} and $(1+t)^{1+2k}$ over $[0,t]$ to achieve
 \begin{align}
\mathcal{E}_k(t) +  \int_0^t  \mathcal{D}_{k}(s) ds
  \le C(\ga, M )  \sum_{0\le \iota \le k} \mathcal{E}_\iota (0).  \label{12-20}
\end{align}

Integrate the product of \ef{3-73} and $\lt(\ln(1+t)\rt)^q(1+t)^{2k}$ over $[0,t]$ to give
\begin{align}
&\lt(\ln(1+t)\rt)^q(1+t)^{2k} \lt(\int_\mathcal{I} r^4  \bar\rho_0 \lt| \pl_t^k  \zeta    \rt|^2 dr + \widehat{\mathfrak{E}}_k(t) \rt) \notag\\
&+\int_0^t \lt(\ln(1+s)\rt)^q(1+s)^{2k} \widehat{\mathfrak{E}}_k(s) ds \notag\\
\le & C(\ga, M) \lt\{ \sum_{0\le \iota\le k}  \mathcal{E}_\iota(0)
+
\lt(\ln(1+t)\rt)^q (1+t)^{-1} \sum_{0\le \iota\le k-1}  \mathcal{E}_\iota(t) \rt. \notag\\
&\lt.+V(t)+\int_0^t \lt(\ln(1+s)\rt)^q \sum_{0\le\iota\le k-1} \mathcal{D}_\iota(s) ds \rt\}\notag\\
\le & C(\ga, M, q) \sum_{0\le \iota\le k}  \mathcal{E}_\iota(0)
+ C(\ga, M) V(t)\le C(\ga, M, q) \sum_{0\le \iota\le k}  \mathcal{E}_\iota(0), \label{12.21}
\end{align}
where
\begin{align*}
V(t)=&\int_0^t \lt(2k\lt(\ln(1+s)\rt)^q + q \lt(\ln(1+s)\rt)^{q-1}   \rt) \lt( (1+s)^{-1}\mathcal{D}_k(s)  \rt. \\
 & \lt. + (1+s)^{2k-1}\int_\mathcal{I} r^4  \bar\rho_0 \lt| \pl_s^k  \zeta    \rt|^2 dr +(1+s)^{-2}\sum_{0\le\iota\le k-1} \mathcal{E}_\iota(s)
\rt)ds.
\end{align*}
Indeed,  the first inequality in \ef{12.21} follows from \ef{12.19b},  and the second from \ef{12-4-j}, and the third from the following estimates. Due to \ef{12-19}, we have
\be\label{12.21-1}
(1+t)^{-1}\mathcal{E}_0(t)\le C(\ga, M)\mathcal{D}_0(t) .
\ee
Clearly, it holds that
\begin{align}
(1+t)^{2j-1}\int_\mathcal{I} r^4 \bar\rho_0    \lt(\pl_t^{j} \zeta \rt)^2    dr= &(1+t)^{2(j-1)+1}\int_\mathcal{I} r^4 \bar\rho_0    \lt(\pl_t^{j-1} \zeta_t \rt)^2    dr\notag\\
\le
&
\mathcal{D}_{j-1}(t), \ \  j=1,2,3,\cdots, \label{12.21-2}
\end{align}
which, together with \ef{12.19}, implies
\begin{align}
(1+t)^{-1}\mathcal{E}_j(t)=&(1+t)^{2j}\widehat{\mathfrak{E}}_j(t)+ (1+t)^{2j-1}\int_\mathcal{I} r^4 \bar\rho_0    \lt(\pl_t^{k} \zeta \rt)^2    dr \notag\\
\le & \mathcal{D}_j(t)+ \mathcal{D}_{j-1}(t), \ \  j=1,2,3, \cdots. \label{12.21-3}
\end{align}
It yields from \ef{12.21-1}-\ef{12.21-3}, \ef{12-4-j} and \ef{12-20} that
\begin{align*}
V(t) \le & \int_0^t \lt(2k\lt(\ln(1+s)\rt)^q + q \lt(\ln(1+s)\rt)^{q-1}   \rt) \\
&\times \lt(C(\ga, M)(1+s)^{-1}\sum_{0\le\iota\le k} \mathcal{D}_\iota(s)+ \mathcal{D}_{k-1}(s)\rt)ds\\
 \le & C(\ga, M, q) \int_0^t \lt( \sum_{0\le\iota\le k} \mathcal{D}_\iota(s)  +  \lt( \lt( \ln(1+s) \rt)^q +1 \rt) \mathcal{D}_{k-1}(s) \rt) ds\\
 \le & C(\ga, M,q )  \sum_{0\le \iota \le k} \mathcal{E}_\iota (0).
\end{align*}
This finishes the proof of \ef{12.21}.

Finally, we integrate the product of  \ef{3-62} and $\lt(\ln(1+t)\rt)^q(1+t)^{1+2k}$ over $[0,t]$ to obtain
\begin{align}
& \lt(\ln(1+t)\rt)^q (1+t)^{1+2k} \widehat{\mathfrak{E}}_k(t)
\notag \\
 & +\int_0^t \lt(\ln(1+s)\rt)^q (1+s)^{1+2k } \int     r^4 \bar\rho_0    \lt(\pl_s^{k} \zeta_s\rt)^2     drds \notag  \\
  \le & C(\ga,M)  \lt\{ \sum_{0\le \iota \le k} \mathcal{E}_\iota (0) +\lt(\ln(1+t)\rt)^q \sum_{0\le \iota\le k-1}  \mathcal{E}_\iota(t) \rt.\notag \\
& +\int_0^t \lt((1+2k)\lt(\ln(1+s)\rt)^q + q \lt(\ln(1+s)\rt)^{q-1}   \rt)  \lt( (1+s)^{2k} \widehat{\mathfrak{E}}_k(s) ds \rt.\notag \\
& \lt. \lt. + \sum_{0\le \iota\le k-1} (1+s)^{-1} \mathcal{E}_\iota(s) \rt) ds  + \int_0^t \lt(\ln(1+s)\rt)^q  \sum_{0\le \iota\le k} (1+s)^{-1} \mathcal{E}_\iota(s) ds \rt\}\notag \\
  \le &  C(\ga,M, q) \lt( \sum_{0\le \iota \le k} \mathcal{E}_\iota (0)
+ \int_0^t \lt( \lt( \lt( \ln(1+s) \rt)^q +1 \rt)  (1+s)^{2k} \widehat{\mathfrak{E}}_k(s) \rt. \rt.\notag \\
 & \lt.  \lt. +\lt( \ln(1+s) \rt)^q (1+s)^{2k-1}\int_\mathcal{I} r^4 \bar\rho_0    \lt(\pl_s^{k} \zeta \rt)^2    dr \rt) ds  \rt)\notag \\
 \le & C(\ga,M, q)   \sum_{0\le \iota \le k} \mathcal{E}_\iota (0) , \label{date3.26}
\end{align}
where the second inequality follows from \ef{12-4-j}, \ef{12.21-1} and \ef{12.21-3}, and the last one from \ef{12.21}, \ef{12.19b}, \ef{12-20}, \ef{12.21-2} and \ef{12-4-j}. This, together with \ef{12.21}, implies that
 \begin{align*}
 \lt( \ln(1+t) \rt)^q  \mathcal{E}_k(t) +  \int_0^t   \lt( \ln(1+s) \rt)^q  \mathcal{D}_{k}(s) ds
  \le C(\ga, M,q )  \sum_{0\le \iota \le k} \mathcal{E}_\iota (0) .
\end{align*}
Using \ef{12-20} and the inequality we just obtained,  we prove that \ef{12-4} holds for $j=k$. This finishes the proof of Lemma \ref{lem1218}.

\subsection{Proof of Theorem \ref{thm2}}
Let $q>0$ be any fixed constant.
It follows from \ef{lem34est-a} and \ef{12.16} that
\begin{align}
  \sum_{j=0,1,2} (1+  t)^{2j} \lt \|\pl_t^j \zeta (\cdot,t) \rt \|_{L^\iy}^2   +  \sum_{j=0,1} (1+  t)^{2j} \lt \|\pl_t^j \pl_r \zeta  (\cdot,t) \rt\|_{L^\iy}^2
 \notag \\
    \le  C(\ga, M, 2q) \lt(\lt(\ln(1+t)\rt)^{2q} +1\rt)^{-1} \mathcal{E}(0)
    \notag \\
    \le  C(\ga, M, q) \lt(\lt(\ln(1+t)\rt)^{q} +1\rt)^{-2} \mathcal{E}(0). \label{March3-3}
\end{align}
Recall that for $t\ge 0$,
\begin{subequations}\label{decayforh}\begin{align}
 0\le h(t) \le C(\ga) (1+t)^{\frac{1}{3\ga-1}-1}\ln(1+t),  \label{March-3}\\
  |h_t(t)| \le C(\ga) (1+t)^{\frac{1}{3\ga-1}-2}\ln(1+t). \label{March-4}
 \end{align}\end{subequations}
This is estimate (2.14) in \cite{HZ}, whose proof follows from the ODE analyses. Note that
\begin{align}
&\rho(\eta(r,t),t)-\bar\rho(\bar\eta(r,t),t)=\frac{r^2 \bar\rho_0(r)}{\eta^2(r,t)\eta_r(r,t)}-\frac{r^2 \bar\rho_0(r)}{\bar\eta^2(r,t)\bar\eta_r(r,t)},\label{March3-1}\\
& u(\eta(r,t),t)-\bar u(\bar\eta(r,t),t)= \eta_t(r,t)- \bar\eta_t(r,t). \label{March3-2}
\end{align}
where \ef{March3-1} follows from \ef{Feb25-1}, \ef{newlagrangiandensity} and \ef{bareta}; and \ef{March3-2} from \ef{Feb26}, \ef{haz2.2} and \ef{bareta}.
Then,  \ef{thm2est1} follows from \ef{bareta}, \ef{March-3} and \ef{March3-3}; and
\ef{thm2est2} from \ef{March-4} and \ef{March3-3}.

For the boundary behavior, \ef{thm2est3} follows from \ef{vbs}, \ef{bareta}, \ef{March-3} and \ef{March3-3}.  \ef{thm2est4} follows from $\tilde\eta=r \tilde\eta_r$, \ef{decayb}, \ef{lem34est-b} and
\bee\label{}\begin{split}
 \frac{d^k R(t)}{dt^k}= \pl_t^{k}\tilde\eta \lt(\sqrt{\frac{A_3}{B_3}},t\rt)  +\lt( r \pl_t^k \zeta\rt)\lt(\sqrt{\frac{A_3}{B_3}},t\rt), \ \ k=1,2,3.
  \end{split}\eee
This finishes the proof of Theorem \ref{thm2}.

\section*{Acknowledgements}

This research was supported in part by NSFC  Grants 11822107 and 11671225. The author would like to thank Professor Chongchun Zeng for helpful discussions.

%\section*{Compliance with Ethical Standards} The author declares that she has no conflict of interest. The author also declares that this research does not involve  human participants or animals, and thus an informed consent is not involved.

\bibliographystyle{plain}

\noindent Huihui Zeng\\
Department of Mathematics\\
\& Yau Mathematical Sciences Center\\
Tsinghua University\\
Beijing, 100084, China;\\
Email: hhzeng@mail.tsinghua.edu.cn

\end{document}